\documentclass[10pt, a4paper, fleqn]{article}
\usepackage[latin1]{inputenc}
\usepackage{amsmath}
\usepackage[numbers]{natbib}
\usepackage{amsfonts}
\usepackage{amssymb}
\usepackage{color}

\allowdisplaybreaks

\numberwithin{equation}{section}
\newtheorem{os}{Remark}[section]
\newtheorem{te}{Theorem}[section]

\author{D'OVIDIO Mirko\footnote{Dipartimento di Statistica, Probabilit\`a e Statistiche Applicate, ''Sapienza'' University of Rome, P.le Aldo Moro n. 5, 00185 Rome (Italy), tel: +390649910499, fax: +39064959241, e-mail:mirko.dovidio@uniroma1.it} \and ORSINGHER Enzo\footnote{Corresponding author: Dipartimento di Statistica, Probabilit\`a e Statistiche Applicate, ''Sapienza'' University of Rome, P.le Aldo Moro n. 5, 00185 Rome (Italy), tel: +390649910585, fax:+39064959241, e-mail:enzo.orsingher@uniroma1.it}}

\title{COMPOSITION OF PROCESSES AND RELATED PARTIAL DIFFERENTIAL EQUATIONS}

\begin{document}

\maketitle

\begin{abstract}
In this paper different types of compositions involving independent fractional Brownian motions $B^j_{H_j}(t)$, $t>0$, $j=1,2$ are examined.\\ 
The partial differential equations governing the distributions of $I_F(t)=B^1_{H_1}(|B^2_{H_2}(t)|)$, $t>0$ and $J_F(t)=B^1_{H_1}(|B^2_{H_2}(t)|^{1/H_1})$, $t>0$ are derived by different methods and compared with those existing in the literature and with those related to $B^1(|B^2_{H_2}(t)|)$, $t>0$. The process of iterated Brownian motion $I^n_F(t)$, $t>0$ is examined in detail and its moments are calculated. Furthermore for $J^{n-1}_F(t)=B^1_{H}(|B^2_H(\ldots |B^n_H(t)|^{1/H} \ldots)|^{1/H})$, $t>0$ the following factorization is proved  $J^{n-1}_F(t)=\prod_{j=1}^{n} B^j_{\frac{H}{n}}(t)$, $t>0$. A series of compositions involving Cauchy processes and fractional Brownian motions are also studied and the corresponding non-homogeneous wave equations are derived.
\end{abstract}


\textbf{AMS:} Primary 60J65, 60J60, 26A33.

\section{Introduction}

In the recent probabilistic literature there are some papers devoted to the interplay between various forms of compositions of different processes and the partial differential equations governing their distributions.\\
The best known example of composition of processes is the iterated Brownian motion (IBM) defined as
\begin{equation}
I(t)=B^{1} \left( |B^{2}(t) | \right), \quad t>0
\label{process:I}
\end{equation}
where $B^1$ and $B^2$ are independent Brownian motions. The IBM has been introduced by Burdzy \cite{BU931} and its properties  like the iterated logarithm law, the fourth-order variation and many others have been analyzed in a series of papers dating back to the middle of Nighties (see Burdzy \cite{BU932}, Khoshnevisan and Lewis \cite{KL96}). The study of iterated Brownian motion has been stimulated by the analysis of diffusions in cracks (see for example, Chudnovsky and Kunin \cite{CK87},  DeBlassie \cite{DB04}). It is well-known that the distribution of \eqref{process:I}, say $q(x,t)$, satisfies the fractional differential equation
\begin{equation}
\frac{\partial^{1/2} q}{\partial t^{1/2}}=\frac{1}{2^{\frac{3}{2}} } \frac{\partial^2 q}{\partial x^2}, \quad x \in \mathbb{R}, \, t >0
\label{e:fractionalOp}
\end{equation}
subject to the initial condition
\begin{equation}
q(x,0)=\delta(x)
\end{equation}
as well as the fourth-order p.d.e. (see DeBlassie \cite{DB04}, Orsingher and Zhao \cite{OZ99})
\begin{equation}
\frac{\partial q}{\partial t}= \frac{1}{2^3} \frac{\partial^4 q}{\partial x^4} + \frac{1}{2 \sqrt{2 \pi t}} \frac{d^2}{d x^2} \delta(x), \quad x \in \mathbb{R}, \, t >0.
\label{eq:allouba}
\end{equation}
Clearly $\frac{d^2}{d x^2}\delta(x)=\delta^{\prime \prime}$ must be understood in the sense that for every test function $\phi$, $<\phi, \delta^{\prime \prime} >=\phi^{\prime \prime}(0)$. This p.d.e. connection is valid for any process $X(|B(t)|)$, $t>0$ where $X$ is a Markov process independent from the Brownian motion $B$ (see Allouba and Zheng \cite{AZ01}, Baeumer et al. \cite{BMN09}, Meerschaert et al. \cite{MNV09}, Baeumer and Meerschaert \cite{BauMee01}). The time-fractional derivative appearing in \eqref{e:fractionalOp} must be understood in the sense 
\begin{equation*}
\frac{\partial^\nu q}{\partial t^\nu}(x,t) = \frac{1}{\Gamma(m-\nu)} \int_{0}^{t} \frac{\frac{\partial^m}{\partial t^m} q(x,s)}{(t-s)^{\nu+1-m}} ds, \quad m-1 < \nu < m. 
\end{equation*}
and
\begin{equation*}
\frac{\partial^\nu q}{\partial t^\nu}(x,t) =\frac{\partial^m q}{\partial t^m}(x,t), \quad \nu=m
\end{equation*}
where $m-1 = \lfloor \nu \rfloor$.\\
For the $n$-times iterated Brownian motion
\begin{equation}
I^n(t)=B^1 (| B^2 (\ldots (|B^{n+1}(t)|) \ldots)|)
\end{equation}
the distribution function solves the fractional equation
\begin{equation}
\frac{\partial^{1/2^n}}{\partial t^{1/2^n}}q=2^{\frac{1}{2^n} - 2} \frac{\partial^2}{\partial x^2}q, \quad x \in \mathbb{R}, \, t >0.
\end{equation}
Extensions of these results to the $d$-dimensional case are treated in Nane \cite{NAa}, Orsingher and Beghin \cite{OB09}. The aim of this paper is to present different types of iterated processes constructed with independent fractional Brownian motions and to explore the partial differential equations governing the corresponding distributions.\\
By considering fractional Brownian motions and Cauchy processes (suitably combined) we obtain qualitatively different equations including a  non-homogeneous wave equation. In the probabilistic literature second-order hyperbolic equations emerge in the study of the telegraph process and of planar random motions with infinite directions. We show here that processes like $C^1(|C^2(t)|)$ ($C^j$, $j=1,2$ are independent Cauchy processes) also are related to wave equations (see also Nane \cite{NA08}). We consider here different types of combinations of processes involving the fractional Brownian motion and the Cauchy process. In particular, we show that for
\begin{equation}
I_{CB_H}(t) = C(|B_H(t)|), \quad t >0
\end{equation}
(where $C$ is a Cauchy process independent from the fractional Brownian motion $B_H$) the probability law is a solution to the non-homogeneous heat equation
\begin{equation}
\frac{\partial q}{\partial t} = \frac{2 H t^{H-1}}{\pi x^2 \sqrt{2 \pi}} - Ht^{2H-1} \frac{\partial^2 q}{\partial x^2}, \quad x \in \mathbb{R}, \, t >0,
\label{eq:densityCB}
\end{equation}
with $H \in (0,1)$. In the special case $H=\frac{1}{2}$ the p.d.e. \eqref{eq:densityCB} becomes the governing equation of the process
\begin{equation}
I_{CB} (t) = C(|B(t)|), \quad t >0.
\end{equation}
For the vector process 
\begin{equation}
I_{BC}^n (t) = \left( B^1(|C(t)|), \ldots , B^n(|C(t)|) \right),  \quad t>0
\end{equation}
we have that the joint law $q=q(x_1, \ldots , x_n, t)$ solves the fourth-order equation 
\begin{equation}
\label{eq:pdeBC}
\frac{\partial^2 q}{\partial t^2} = -\frac{1}{2^2} \triangle^2 q -\frac{1}{\pi t} \triangle \delta(x_1) \ldots \delta(x_n), \quad x \in \mathbb{R}^n, \, t >0. 
\end{equation}
This fact has been noted by Nane in \cite{NA08}, \cite{NAa}. It seems that it is difficult to obtain the equation corresponding to
\begin{equation}
\left( B^1_{H_1}(|C(t)|), \ldots , B^n_{H_n}(|C(t)|) \right), \quad t>0
\end{equation}
where the independent fractional Brownian motions $B^j_{H_j}(t)$ are involved.\\
Finally we prove that the composition of two independent Cauchy processes 
\begin{equation}
I_{C}(t)=C_1(|C_2(t)|), \quad t>0
\label{process:CC}
\end{equation}
has a probability law satisfying the non-homogeneous wave equation
\begin{equation}
\label{eq:pdeCC2}
\frac{\partial^2 q}{\partial t^2} = -\frac{2}{\pi^2 t x^2} + \frac{\partial^2 q}{\partial x^2}, \quad x \in \mathbb{R}, \, t >0
\end{equation}
which, in our view, is the most striking result of all these combinations of well-known processes.\\
We also prove that the one-dimensional distribution of \eqref{process:CC} coincides with the multiplicative process 
\begin{equation}
C_1(|C_2(t)|) \stackrel{i.d.}{=} \frac{1}{2} C_1(\sqrt{2t})C_2(\sqrt{2t}), \quad t >0.
\label{processCCdist}
\end{equation}
It is also true that the composition of Cauchy processes satisfies the following curious relationship
\begin{equation}
C^1(|C^2(t)|) \stackrel{i.d.}{=} \frac{1}{C^1\left(|C^2\left(\frac{1}{t}\right)|\right)}, \quad t>0
\end{equation}
We remark that if the guiding process is a Cauchy process the related equation has a forcing function defined on the whole space-time domain, while in the case of a guiding process coinciding with Brownian motion, the forcing term is concentrated on the $x$ line.\\
The core of our paper concerns the iterated fractional Brownian motion whose general form is given by
\begin{equation}
I_F(t)=B^{1}_{H_1} (|B^{2}_{H_2}(t)|), \quad t>0
\label{process:IF}
\end{equation}
where the processes $B^j_{H_j}(t)$, $t>0$ have covariance function 
\begin{equation}
E \left\lbrace B^{j}_{H_j} (s) B^{j}_{H_j} (t) \right\rbrace = \frac{1}{2} \left( |t|^{2H_j} + |s|^{2H_j} - |t-s|^{2H_j} \right), \quad j=1,2
\end{equation}
and $0<H_1, H_2 < 1$. \\

For the iterated fractional Brownian motion 
\begin{equation}
\tilde{I}_H(t)=B^1(|B^2_H(t)|), \quad t>0, \; H \in (0,1)
\label{process:Itilde}
\end{equation}
where $B^1$ is a standard Brownian motion independent from the fractional Brownian motion $B^2_H$, we show that the equation governing its distribution $q=q(x,t)$ is
\begin{equation}
\frac{\partial q}{\partial t} = \frac{H t^{2H - 1}}{2^2} \frac{\partial^4 q}{\partial x^4} + \frac{H t^{2H - 1}}{2 \sqrt{2 \pi t^{2H}}} \frac{d^2}{d x^2} \delta(x), \quad x \in \mathbb{R}, \, t >0, \quad H \in (0,1)
\label{eq:IFF}
\end{equation}
with $q(x,0)=\delta(x)$.\\
The introduction of the fractional Brownian motion makes the equation \eqref{eq:IFF} slightly different from \eqref{eq:allouba} in that it has non-constant coefficients. In the case where \eqref{process:Itilde} is replaced by the related process
\begin{equation}
L_F(t)=t^K B^1(|B^2_H(t)|) \quad  t>0, \, K>0, \, H \in (0,1)
\end{equation}
the corresponding equation becomes a bit different and has the form
\begin{equation}
t\frac{\partial q}{\partial t} = -K \frac{\partial}{\partial x} \left( xq \right) + \frac{H}{4} t^{4K + 2H} \frac{\partial^4 q}{\partial x^4} + \frac{Ht^{2K +H}}{\sqrt{2\pi}} \frac{d^2}{d x^2} \delta(x), \quad x \in \mathbb{R}, \, t>0, \, K \geq 0
\end{equation}
with an additional lower-order term with space-dependent coefficients.\\
For the general iterated fractional Brownian motion 
\begin{equation}
I_F(t)=B^1_{H_1}(|B^2_{H_2}(t)|), \quad t>0, \; 0 < H_1,H_2 < 1
\label{process:IFH}
\end{equation}
we have arrived at a substantially different equation which shows the fundamental role of the process $B_{H_1}$. The equation pertaining to process \eqref{process:IFH} reads
\begin{equation}
(1+H_1 H_2) t \frac{\partial p}{\partial t} + t^2 \frac{\partial^2 p}{\partial t^2} =H_1^2 H_2^2 \left\lbrace 2x\frac{\partial p}{\partial x} + x^2 \frac{\partial^2 p}{\partial x^2} \right\rbrace, \quad x \in \mathbb{R}, \, t >0.
\label{eqORA}
\end{equation}
We note that equation \eqref{eqORA} has almost the same structure of (5.2) of Orsingher and De Gregorio \cite{ODG07} which emerges in the study of random motions at finite velocity in hyperbolic spaces.\\
For the process \eqref{process:IFH} we have the general expression of even-order moments
\begin{equation}
E\left\lbrace B^1_{H_1}(|B^2_{H_2}(t)|) \right\rbrace^{2k} = 2^{2-k-kH_1} t^{2kH_1H_2} \frac{\Gamma(2k)}{\Gamma(k)} \frac{\Gamma(2kH_1)}{\Gamma(kH_1)}
\end{equation}
which for $H_1=H_2=1/2$ yields
\begin{equation}
E\left\lbrace B^1(|B^2(t)|) \right\rbrace^{2k} = \frac{2^{k/2}}{2^{2k}} \frac{(2k)!}{\Gamma\left(\frac{k}{2} +1 \right)} t^{k/2}
\end{equation} 
which is the first formula of the remark 3.3 of Orsingher and Beghin \cite{OB09}.\\
Similarly for the $n$-times iterated fractional Brownian motion we have that
\begin{equation}
E\left\lbrace B^1_{H_1}(|B^2_{H_2}(\ldots |B^{n+1}_{H_{n+1}}(t)| \ldots)|) \right\rbrace^{2k} =  \frac{2^{n+1} t^{2k\prod_{j=1}^{n+1}H_j}}{2^{k\sum_{r=0}^{n} \prod_{j=0}^{r}H_j}} \prod_{r=0}^{n} \frac{\Gamma(2k\prod_{j=0}^{r}H_j)}{\Gamma(k\prod_{j=0}^{r}H_j)}
\end{equation}
for $n \geq 1$, $k \in \mathbb{N}$, $H_0=1$ and $H_j \in (0,1)$ for $j=1,2 \ldots ,n$ .\\ 
Another process obtained by composing independent fractional Brownian motions is
\begin{equation}
J_F^1 (t) = B_H^{1} (| B_H^{2}(t) |^\frac{1}{H}), \quad t>0, \; H \in (0,1)
\label{process:IIF}
\end{equation}
which has proved to be much more tractable than \eqref{process:IF}. For \eqref{process:IIF} it is even possible to obtain the explicit law in the following form
\begin{align}
q(x,t) = & 2 \int_{0}^{\infty} \frac{e^{-\frac{x^2}{2 s^2}}}{\sqrt{2 \pi s^2}} \frac{e^{-\frac{s^2}{2 t^{2H}}}}{\sqrt{2 \pi t^{2H}}} ds = \frac{1}{\pi t^H} K_0 \left( \frac{|x|}{t^H} \right), \quad H \in (0,1], \; x \in \mathbb{R} \setminus \{ 0 \}
\label{eq:densityIIF}
\end{align}
where $K_0(x)$ is the Bessel modified  function defined as
\begin{equation}
K_0(x)=\int_{0}^{\infty} \frac{1}{s} e^{-\frac{x^2}{4 s^2} - s^2} ds, \quad x \neq 0.
\end{equation}
The p.d.e. governing the distribution (\ref{eq:densityIIF}) is 
\begin{equation}
\label{eq:problemK0}
\frac{\partial}{\partial t} q(x,t)=-Ht^{2H-1} \left( 2\frac{\partial^2}{\partial x^2} + x \frac{\partial^3}{\partial x^3} \right) q(x,t)
\end{equation}
\begin{equation*}
\textcolor{white}{\frac{\partial}{\partial t} q(x,t)}= -Ht^{2H-1} \frac{\partial^2}{\partial x^2} \left( x\frac{\partial}{\partial x} q(x,t) \right), \quad x \in \mathbb{R}, \, t >0.
\end{equation*}
We have also considered an extension of \eqref{process:IIF} by taking into account the process
\begin{equation}
J^2_F (t) = B_H^{1} (| B_H^{2}(| B_H^{3}(t) |^{\frac{1}{H}}) |^\frac{1}{H}), \quad t>0, \; H \in (0,1)
\end{equation}
whose distribution is related to a fourth-order p.d.e. with non-constant coefficients which reads
\begin{equation}
\frac{\partial q}{\partial t} = H t^{2H-1} \left\lbrace 4 \frac{\partial^2 }{\partial x^2} + 5x \frac{\partial^3 }{\partial x^3} + x^2 \frac{\partial^4 }{\partial x^4} \right\rbrace q, \quad x \in \mathbb{R}, \, t >0.
\end{equation}
We note that for the extension of \eqref{process:IIF}, say $J_F^n(t)$, the following equivalences in distribution hold
\begin{equation} 
\label{JnDistribution1}
J^n_F (t)= B_H^{1} (| B_H^{2}( \ldots | B^{n+1}_{H}(t) |^{\frac{1}{H}} \ldots ) |^\frac{1}{H}) \stackrel{i.d.}{=} B^{1} (| B^{2}( \ldots | B^{n+1}_{H}(t) |^2 \ldots ) |^2), \quad t>0.
\end{equation}
We will also prove that for $J_F^{n-1}(t)$, $t>0$ the factorization
\begin{equation}
J^{n-1}_F (t) \stackrel{i.d.}{=} \prod_{i=1}^{n} B^i_{\frac{H}{n}} (t), \quad t>0, \; H \in (0,1) 
\label{prodBH}
\end{equation}
in terms of independent fractional Brownian motions $B^i_{\frac{H}{n}}$, $i=1,2 \ldots , n$ holds.
In view of decomposition \eqref{prodBH} one can also observe that 
\begin{equation} 
B_H^{1} (| B_H^{2}(t) |^\frac{1}{H}) \stackrel{i.d.}{=} B^1(| B_H^{2}(t) |^2) \stackrel{i.d.}{=} B^1_{\frac{H}{2}}(t) B^2_{\frac{H}{2}}(t), \quad t>0, \, H \in (0,1].
\label{JnDistribution2}
\end{equation}
We have considered some multidimensional compositions of Brownian motions. In Orsingher and Beghin \cite{OB09}, Nane \cite{NAa} it is proved that the $n$-dimensional vector process 
\begin{equation} 
\left( B^1(|B(t)|), \ldots , B^n(|B(t)|) \right), \quad t>0
\label{eqBORA}
\end{equation}
(where $B^1, \ldots , B^n, B$ are independent Brownian motions) has a distribution function which is a solution to
\begin{equation}
\frac{\partial^{1/2} u}{\partial t^{1/2}} = \frac{1}{2^{3/2}} \sum_{j=1}^{n} \frac{\partial^2 u}{\partial x^2_j}, \quad u(x_1, \ldots , x_n, 0)=\prod_{j=1}^{n} \delta(x_j), \quad x \in \mathbb{R}^n, \, t >0.
\end{equation}
In Allouba and Zheng \cite{AZ01}, DeBlassie \cite{DB04} it is shown that the law of \eqref{eqBORA} satisfies the fourth-order equation
\begin{equation}
\frac{\partial u}{\partial t} = \frac{1}{2^3} \left( \sum_{j=1}^{n} \frac{\partial^2}{\partial x^2_j} \right)^2 u + \frac{1}{2\sqrt{2\pi t}} \sum_{j=1}^{n} \frac{d^2}{d x^2_j} \prod_{r=1}^{n} \delta(x_r), \quad x \in \mathbb{R}^n, \, t>0.
\label{BORA2}
\end{equation}
For the $n$-dimensional process $(B^1(|B_{H}(t)|), \ldots , B^n(|B_{H}(t)|))$ we arrive at the fourth-order p.d.e.
\begin{equation}
\frac{\partial q}{\partial t}=\frac{H t^{2H-1}}{\sqrt{2 \pi t}} \triangle \prod_{j=1}^{n}\delta(x_j) + \frac{H t^{2H-1}}{2^2}\triangle^2 q, \quad x \in \mathbb{R}^n, \, t>0
\end{equation}
which includes \eqref{BORA2} as a special case for $H=\frac{1}{2}$. 

\section{The iterated fractional Brownian motion}

For the iterated Brownian motion $I(t)=B^1(|B^2(t)|)$, $t >0$ it has been shown by different authors that the distribution
\begin{equation}
p(x,t)=2 \int_{0}^{\infty} \frac{e^{-\frac{x^2}{2s}}}{\sqrt{2 \pi s}} \frac{e^{-\frac{s^2}{2t}}}{\sqrt{2 \pi t}} ds, \quad x \in \mathbb{R}, \; t>0
\label{eq:dist1}
\end{equation}
with initial condition $p(x,0)=\delta(x)$ is the solution to the fourth-order equation
\begin{equation}
\frac{\partial p}{\partial t} = \frac{1}{2^3} \frac{\partial^4 p}{\partial x^4} + \frac{1}{2 \sqrt{2 \pi t}} \frac{d^2}{d x^2} \delta(x), \quad x \in \mathbb{R}, \, t >0
\label{eq:pde1}
\end{equation}
(see for example Allouba and Zheng \cite{AZ01}, DeBlassie \cite{DB04}, Nane \cite{NA08}).\\
The distribution (\ref{eq:dist1}) is also the solution to the time-fractional diffusion equation
\begin{equation}
\frac{\partial^{1/2} q}{\partial t^{1/2}}= \frac{1}{2^{3/2}} \frac{\partial^2 q}{\partial x^2}, \quad x \in \mathbb{R}, \, t>0
\end{equation}
subject to the initial condition $q(x,0)=\delta(x)$ (see Orsingher and Beghin \cite{OB04}, \cite{OB09}).\\
We now show below that the law
\begin{equation}
p(x,t)=2 \int_{0}^{\infty} \frac{e^{-\frac{x^2}{2s^{2H_1}}}}{\sqrt{2 \pi s^{2H_1}}} \frac{e^{-\frac{s^2}{2t^{2H_2}}}}{\sqrt{2 \pi t^{2H_2}}} ds, \quad x \in \mathbb{R}, \; t>0, \; H_1, H_2 \in (0,1)
\end{equation}
the iterated fractional Brownian motion $I_F(t)=B^1_{H_1}(|B^2_{H_2}(t)|)$, $t>0$ is a solution to a second-order differential equation with non-constant coefficients.\\
We first need the following auxiliary results concerning the density of fractional Brownian motion.
\begin{te}
\label{te:BH}
The fractional Brownian motion $B_{H}(t)$, $t >0$, $H \in (0,1)$ has a probability density $p_{H}(x,t)$, $x \in \mathbb{R}$, $t>0$ which solves the heat equation
\begin{equation}
\label{eq:problemBH}
\frac{\partial}{\partial t}p_{H}=Ht^{2H-1}\frac{\partial^2}{\partial x^2}p_{H}
\end{equation}
with initial condition $p_H(x,0)=\delta(x)$.
\end{te}
\paragraph{Proof:}
The distribution of $B_{H}(t)$, $t>0$ can be written as inverse Fourier transform
\begin{equation}
p_{H}(x,t)=\frac{1}{2\pi} \int_{-\infty}^{\infty} e^{-i\beta x} e^{-\frac{1}{2}\beta^2 t^{2H}} d\beta.
\end{equation}
A simple calculation proves that (\ref{eq:problemBH}) holds. For $H=\frac{1}{2}$ equation (\ref{eq:problemBH}) coincides with the classical heat-equation. 
\begin{flushright}
$\blacksquare$
\end{flushright}

\begin{os}
\normalfont
For the Gaussian law
\begin{equation}
q(x,t)=\frac{e^{-\frac{x^2}{2g(t)}}}{\sqrt{2\pi g(t)}}, \quad x \in \mathbb{R}, \; t>0
\label{eq:densityB1}
\end{equation}
with $g \in C^1$, $g:(0,\infty) \mapsto (0,\infty)$, we can easily accertain that the governing equation is
\begin{equation}
\label{eq:problemB1}
\frac{\partial}{\partial t}q= \frac{dg}{dt} \frac{1}{2} \frac{\partial^2}{\partial x^2}q, \quad x \in \mathbb{R}, \; t>0 .
\end{equation}
\end{os}
 
We prove below that the processes $\tilde{I}_H(t)=B^1(|B^2_H(t)|)$ and its slightly extended version $L_H(t)=t^K B^1(|B^2_H(t)|)$ have a distribution satisfying a fourth-order p.d.e. somewhat similar to \eqref{eq:pde1}. We analyze in the next theorem the law of $L_H(t)$ and extract the related equation for $\tilde{I}_H(t)$ by assuming that $K=0$.
\begin{te}
The law of the process 
\begin{equation}
L_F(t)=t^K B^1(|B^2_H(t)|), \quad  t>0, \, K>0
\end{equation}
is solution of the fourth-order equation
\begin{equation}
t\frac{\partial q}{\partial t}=-K\frac{\partial}{\partial x}\left( xq \right) +\frac{H}{4} t^{4K+2H} \frac{\partial^4 q}{\partial x^4} + \frac{H t^{2K+H}}{\sqrt{2\pi}} \frac{\partial^2}{\partial x^2}\delta(x), \quad x \in \mathbb{R}, \, t>0, \, K>0
\end{equation}
subject to the initial condition $q(x,0)=\delta(x)$.
\end{te}
\paragraph{Proof:}
We start by writing down the distribution of $L_F(t)$, $t>0$ as
\begin{equation}
p(x,t)=2 \int_{0}^{\infty} \frac{e^{-\frac{x^2}{2s t^{2K}}}}{\sqrt{2 \pi s t^{2K}}} \frac{e^{-\frac{s^2}{2 t^{2H}}}}{\sqrt{2 \pi t^{2H}}} ds, \quad x \in \mathbb{R}, \; t>0, \; K \geq 0.
\label{eq:densityphi}
\end{equation}
By taking the time derivative of \eqref{eq:densityphi} we have that
\begin{equation}
\frac{\partial p}{\partial t} =2 \int_{0}^{\infty} \frac{\partial}{\partial t} \left( \frac{e^{-\frac{x^2}{2 s t^{2K}}}}{\sqrt{2 \pi s t^{2K}}} \right) \frac{e^{-\frac{s^2}{2 t^{2H}}}}{\sqrt{2 \pi t^{2H}}} ds + 2 \int_{0}^{\infty} \frac{e^{-\frac{x^2}{2s t^{2K}}}}{\sqrt{2 \pi s t^{2K}}} \frac{\partial}{\partial t} \left( \frac{e^{-\frac{s^2}{2 t^{2H}}}}{\sqrt{2 \pi t^{2H}}} \right) ds
\label{dimLH}
\end{equation}
\begin{equation*}
\textcolor{white}{\frac{\partial p}{\partial t}}= \int_{0}^{\infty} \frac{\partial^2}{\partial x^2} \left( \frac{e^{-\frac{x^2}{2 s t^{2K}}}}{\sqrt{2 \pi s t^{2K}}} \right) 2K t^{2K-1} s \frac{e^{-\frac{s^2}{2 t^{2H}}}}{\sqrt{2 \pi t^{2H}}} ds 
\end{equation*}
\begin{equation*}
\textcolor{white}{\frac{\partial p}{\partial t}}+ 2H t^{2H-1} \int_{0}^{\infty} \frac{e^{-\frac{x^2}{2s t^{2K}}}}{\sqrt{2 \pi s t^{2K}}} \frac{\partial^2}{\partial s^2} \left( \frac{e^{-\frac{s^2}{2 t^{2H}}}}{\sqrt{2 \pi t^{2H}}} \right) ds.
\end{equation*}
In the last step we used equations  \eqref{eq:problemBH} and \eqref{eq:problemB1}. The first integral in \eqref{dimLH} can be developed as
\begin{equation*}
-\frac{K}{t} \frac{\partial}{\partial x} \left\lbrace x\, 2 \int_{0}^{\infty} \frac{e^{-\frac{x^2}{2 s t^{2K}}}}{\sqrt{2 \pi s t^{2K}}} \frac{e^{-\frac{s^2}{2t^{2H}}}}{\sqrt{2 \pi t^{2H}}} ds \right\rbrace  = -\frac{K}{t} \frac{\partial}{\partial x} \left( xp \right).
\end{equation*}
By integrating by parts the second integral in \eqref{dimLH} we have that
\begin{equation*}
2H t^{2H -1} \frac{e^{-\frac{x^2}{2 s t^{2K}}}}{\sqrt{2 s \pi t^{2K}}} \frac{\partial}{\partial s}\left( \frac{e^{-\frac{s^2}{2 t^{2H}}}}{\sqrt{2 \pi t^{2H}}} \right) \Bigg|_{s=0}^{s=\infty} 
\end{equation*}
\begin{equation*}
- 2Ht^{2H-1} \int_{0}^{\infty} \frac{\partial}{\partial s} \left( \frac{e^{-\frac{x^2}{2 s t^{2K}}}}{\sqrt{2 \pi s t^{2K}}} \right) \frac{\partial}{\partial s} \left( \frac{e^{-\frac{s^2}{2 t^{2H}}}}{\sqrt{2 \pi t^{2H}}} \right) ds =
\end{equation*}
\begin{equation*}
-2H t^{2H -1} \frac{\partial}{\partial s} \left( \frac{e^{-\frac{x^2}{2 s t^{2K}}}}{\sqrt{2 s \pi t^{2K}}} \right) \frac{e^{-\frac{s^2}{2 t^{2H}}}}{\sqrt{2 \pi t^{2H}}} \Bigg|_{s=0}^{s=\infty} +
\end{equation*}
\begin{equation*}
2Ht^{2H-1} \int_{0}^{\infty} \frac{\partial^2}{\partial s^2} \left( \frac{e^{-\frac{x^2}{2 s t^{2K}}}}{\sqrt{2 \pi s t^{2K}}} \right) \frac{e^{-\frac{s^2}{2 t^{2H}}}}{\sqrt{2 \pi t^{2H}}} ds.
\end{equation*}
By assuming $\phi=st^{2K}$ and in view of \eqref{eq:problemB1} we can write that
\begin{equation*}
\frac{\partial}{\partial s} \frac{e^{-\frac{x^2}{2 s t^{2K}}}}{\sqrt{2 \pi s t^{2K}}} = \frac{t^{2K}}{2}\frac{\partial^2}{\partial x^2} \frac{e^{-\frac{x^2}{2 s t^{2K}}}}{\sqrt{2 \pi s t^{2K}}} \quad \textrm{and} \quad \frac{\partial^2}{\partial s^2} \frac{e^{-\frac{x^2}{2 s t^{2K}}}}{\sqrt{2 \pi s t^{2K}}} = \frac{t^{4K}}{2^2} \frac{\partial^4}{\partial x^4} \frac{e^{-\frac{x^2}{2 s t^{2K}}}}{\sqrt{2 \pi s t^{2K}}}. 
\end{equation*}
All this leads to
\begin{equation*}
\frac{\partial p}{\partial t} = - \frac{K}{t} \frac{\partial}{\partial x} \left( xp \right)
\end{equation*}
\begin{equation*}
\textcolor{white}{\frac{\partial p}{\partial t} =} - Ht^{2H + 2K -1} \frac{\partial^2}{\partial x^2} \left( \frac{e^{-\frac{x^2}{2 s t^{2K}}}}{\sqrt{2 \pi s t^{2K}}} \right) \frac{e^{-\frac{s^2}{2 t^{2H}}}}{\sqrt{2 \pi t^{2H}}} \Bigg|_{0}^{\infty}
\end{equation*}
\begin{equation*}
\textcolor{white}{\frac{\partial p}{\partial t} =} + \frac{H}{2^2} t^{4K+2H-1} \frac{\partial^4}{\partial x^4} \left\lbrace 2 \int_{0}^{\infty} \frac{e^{-\frac{x^2}{2s t^{2K}}}}{\sqrt{2 \pi s t^{2K}}} \frac{e^{-\frac{s^2}{2 t^{2H}}}}{\sqrt{2 \pi t^{2H}}} ds \right\rbrace 
\end{equation*}
\begin{equation*}
\textcolor{white}{\frac{\partial p}{\partial t}} =-\frac{K}{t}\frac{\partial}{\partial x} \left(xp \right) + \frac{H t^{H+2K-1}}{\sqrt{2\pi}} \frac{d^2}{d x^2} \delta(x) + \frac{H}{2^2} t^{4K+2H-1} \frac{\partial^4 p}{\partial x^4}.
\end{equation*}
This concludes our proof.
\begin{flushright}
$\blacksquare$
\end{flushright} 

\begin{os}
\normalfont
For $K=0$, equation \eqref{eq:densityphi} simplifies and becomes
\begin{equation}
t \frac{\partial q}{\partial t} = \frac{H}{2^2} t^{2H}\frac{\partial^4 q}{\partial x^4} + \frac{H t^{H}}{\sqrt{2 \pi}} \frac{d^2}{d x^2}\delta(x), \quad x \in \mathbb{R}, \, t >0\; H \in (0,1).
\label{pde2}
\end{equation}
Clearly for $H=1/2$ we reobtain DeBlassie's result
\begin{equation}
\frac{\partial q}{\partial t} = \frac{1}{2^3} \frac{\partial^4 q}{\partial x^4} + \frac{ 1}{2\sqrt{2 \pi t}} \frac{d^2}{d x^2}\delta(x), \quad x \in \mathbb{R}, \, t >0.
\end{equation}

It is useful to have a look at table 1, for comparing the governing equations examined so far.

\begin{table}
\centering
\begin{tabular}{c|c|c|c}
Process & Governing & Iterated & Governing \\
& equation & process & equation \\
\hline\hline
$B(t)$ & $\frac{\partial q}{\partial t}=\frac{1}{2}\frac{\partial^2 q}{\partial x^2}$ & $B(|B(t)|)$ & $\frac{\partial q}{\partial t}=\frac{1}{2^3}\frac{\partial^4 q}{\partial x^4} + \frac{1}{\sqrt{2\pi t}} \frac{d^2}{d x^2} \delta$  \\[1ex]
 &  &  & $\frac{\partial^{1/2} q}{\partial t^{1/2}}=\frac{1}{2^{3/2}}\frac{\partial^2 q}{\partial x^2}$ \\[1.5ex]
\hline 
$B_H(t)$ & $\frac{\partial q}{\partial t}=H t^{2H-1}\frac{\partial^2 q}{\partial x^2}$ & $t^K B(|B_H(t)|)$ & $t\frac{\partial q}{\partial t} = -K \frac{\partial}{\partial x} \left( xq \right) + \frac{H}{4} t^{4K + 2H} \frac{\partial^4 q}{\partial x^4}$\\[1ex]
 &  &  & $+ \frac{Ht^{2K +H}}{\sqrt{2\pi}} \frac{d^2}{d x^2} \delta(x)$\\
\end{tabular}
\begin{center}
\textbf{Table 1}
\end{center}
\end{table}
\end{os}
It is easy to check that the vector process $(B^1_{H_1}(t), \ldots B^n_{H_n}(t))$, $t>0$ has a law which satisfies the heat equation
\begin{equation}
\frac{\partial q}{\partial t} = \frac{1}{2} \sum_{j=1}^{n} H_j t^{2H_j -1} \frac{\partial^2 q}{\partial x^2_j}.
\end{equation}
For the $n$-dimensional vector process $(B^1(|B_{H}(t)|), \ldots B^n(|B_{H}(t)|))$ the distribution function can be written as
\begin{equation}
q(x_1, \ldots , x_n, t)=2 \int_{0}^{\infty} ds \int_{\mathbb{R}^n} e^{-i \sum_{j=1}^{n} \alpha_j x_j - \frac{s}{2} \sum_{j=1}^{n} \alpha^2_j} \frac{e^{-\frac{s^2}{2t^{2H}}}}{\sqrt{2 \pi t^{2H}}} d\alpha_1 \ldots d\alpha_n.
\end{equation}
By taking the time derivative and considering \eqref{eq:problemBH} we have that
\begin{equation*}
\frac{\partial q}{\partial t} = 2 \int_{0}^{\infty} \frac{ds}{(2\pi)^n} \int_{\mathbb{R}^n} e^{-i \sum_{j=1}^{n} \alpha_j x_j - \frac{s}{2} \sum_{j=1}^{n} \alpha^2_j} H t^{2H-1} \frac{\partial^2}{\partial s^2} \left( \frac{e^{-\frac{s^2}{2t^{2H}}}}{\sqrt{2 \pi t^{2H}}} \right) d\alpha_1 \ldots d\alpha_n
\end{equation*}
\begin{equation*}
\textcolor{white}{\frac{\partial q}{\partial t}} = \frac{2}{(2\pi)^n} \int_{\mathbb{R}^n} e^{-i \sum_{j=1}^{n} \alpha_j x_j - \frac{s}{2} \sum_{j=1}^{n} \alpha^2_j} H t^{2H-1} \frac{\partial}{\partial s} \left( \frac{e^{-\frac{s^2}{2t^{2H}}}}{\sqrt{2 \pi t^{2H}}} \right) d\alpha_1 \ldots d\alpha_n \Bigg|_{s=0}^{s=\infty}
\end{equation*}
\begin{equation*}
\textcolor{white}{\frac{\partial q}{\partial t}} +  H t^{2H-1} \int_{0}^{\infty} \frac{ds}{(2\pi)^n} \int_{\mathbb{R}^n}  \sum_{j=1}^{n} \alpha^2_j e^{-i \sum_{j=1}^{n} \alpha_j x_j - \frac{s}{2} \sum_{j=1}^{n} \alpha^2_j} \frac{\partial}{\partial s} \left( \frac{e^{-\frac{s^2}{2t^{2H}}}}{\sqrt{2 \pi t^{2H}}} \right) d\alpha_1 \ldots d\alpha_n
\end{equation*}
\begin{equation*}
\textcolor{white}{\frac{\partial q}{\partial t}} =  H t^{2H-1} \int_{0}^{\infty} \frac{ds}{(2\pi)^n} \int_{\mathbb{R}^n} \sum_{j=1}^{n} \alpha^2_je^{-i \sum_{j=1}^{n} \alpha_j x_j - \frac{s}{2} \sum_{j=1}^{n} \alpha^2_j} \frac{\partial}{\partial s} \left( \frac{e^{-\frac{s^2}{2t^{2H}}}}{\sqrt{2 \pi t^{2H}}} \right) d\alpha_1 \ldots d\alpha_n
\end{equation*}
\begin{equation*}
\textcolor{white}{\frac{\partial q}{\partial t}} =  H t^{2H-1} \frac{1}{(2\pi)^n} \int_{\mathbb{R}^n} \sum_{j=1}^{n} \alpha^2_j e^{-i \sum_{j=1}^{n} \alpha_j x_j - \frac{s}{2} \sum_{j=1}^{n} \alpha^2_j} \frac{e^{-\frac{s^2}{2t^{2H}}}}{\sqrt{2 \pi t^{2H}}}d\alpha_1 \ldots d\alpha_n \Bigg|_{s=0}^{s=\infty}
\end{equation*}
\begin{equation*}
\textcolor{white}{\frac{\partial q}{\partial t}} + H t^{2H-1} \int_{0}^{\infty} \frac{ds}{(2\pi)^n} \int_{\mathbb{R}^n} \left( \sum_{j=1}^{n} \alpha^2_j \right)^2 e^{-i \sum_{j=1}^{n} \alpha_j x_j - \frac{s}{2} \sum_{j=1}^{n} \alpha^2_j} \frac{e^{-\frac{s^2}{2t^{2H}}}}{\sqrt{2 \pi t^{2H}}} d\alpha_1 \ldots d\alpha_n
\end{equation*}
\begin{equation*}
\textcolor{white}{\frac{\partial q}{\partial t}} = \frac{Ht^{2H-1}}{\sqrt{2\pi t^{2H}}} \sum_{j=1}^{n} \frac{d^2}{d x^2_j} \delta(x_j) + \frac{Ht^{2H-1}}{2^2} \left( \sum_{j=1}^{n} \frac{\partial^2}{\partial x^2_j} \right) \left( \sum_{j=1}^{n} \frac{\partial^2}{\partial x^2_j}  \right) q(x_1, \ldots , x_n, t)
\end{equation*}
\begin{equation*}
\textcolor{white}{\frac{\partial q}{\partial t}} = \frac{Ht^{2H-1}}{\sqrt{2\pi t^{2H}}} \triangle \delta(x_1) \cdots \delta(x_n) + \frac{Ht^{2H-1}}{2^2} \triangle^2 q
\end{equation*}
where $\triangle=\sum_{j=1}^{n} \frac{\partial^2}{\partial x^2_j}$.\\

The iterated fractional Brownian motion $I^1_F(t)=B^1_{H_1}(|B^2_{H_2}(t)|)$ brings about substantially different p.d.e.'s and we start our analysis from the first-order equation of the next theorem. For a duscussion on the uniqueness of Cauchy problems related to the above equations see Baeumer et al. \cite{BMN09}.
\begin{te}
The distribution of the iterated fractional Brownian motions $I^1_F(t)=B_{H_1}^{1}(|B_{H_2}^{2}(t)|)$ is a solution of the first-order p.d.e.
\begin{equation}
t\frac{\partial p}{\partial t} =- H_1 H_2 \frac{\partial}{\partial x} \left(x p \right), \qquad x \in \mathbb{R}, \; t>0
\label{eq:firstOrd}
\end{equation}
where $H_1$, $H_2 \in (0,1)$.
\end{te}
\paragraph{Proof:} The time derivative of
\begin{equation}
p(x,t)=2 \int_{0}^{\infty} \frac{e^{-\frac{x^2}{2s^{2H_1}}}}{\sqrt{2\pi s^{2H_1}}} \frac{e^{-\frac{s^2}{2t^{2H_2}}}}{\sqrt{2\pi t^{2H_2}}}ds
\end{equation}
writes
\begin{equation*}
\frac{\partial}{\partial t}p(x,t)=2 \int_{0}^{\infty} \frac{e^{-\frac{x^2}{2s^{2H_1}}}}{\sqrt{2\pi s^{2H_1}}} \frac{\partial}{\partial t} \frac{e^{-\frac{s^2}{2t^{2H_2}}}}{\sqrt{2\pi t^{2H_2}}}ds
\end{equation*}
\begin{equation*}
\textcolor{white}{\frac{\partial}{\partial t}p(x,t)}=2H_2 t^{2H_2 -1} \int_{0}^{\infty} \frac{e^{-\frac{x^2}{2s^{2H_1}}}}{\sqrt{2\pi s^{2H_1}}} \frac{\partial^2}{\partial s^2} \frac{e^{-\frac{s^2}{2t^{2H_2}}}}{\sqrt{2\pi t^{2H_2}}}ds
\end{equation*}
\begin{equation*}
\textcolor{white}{\frac{\partial}{\partial t}p(x,t)}=2H_2 t^{2H_2 -1} \left\lbrace \frac{e^{-\frac{x^2}{2s^{2H_1}}}}{\sqrt{2\pi s^{2H_1}}} \frac{\partial}{\partial s} \frac{e^{-\frac{s^2}{2t^{2H_2}}}}{\sqrt{2\pi t^{2H_2}}} \Big|_{s=0}^{s=\infty} - \int_{0}^{\infty} \frac{\partial}{\partial s} \frac{e^{-\frac{x^2}{2s^{2H_1}}}}{\sqrt{2\pi s^{2H_1}}} \frac{\partial}{\partial s} \frac{e^{-\frac{s^2}{2t^{2H_2}}}}{\sqrt{2\pi t^{2H_2}}}ds \right\rbrace 
\end{equation*}
\begin{equation*}
\textcolor{white}{\frac{\partial}{\partial t}p(x,t)}=-2H_2 t^{2H_2 -1}  \int_{0}^{\infty} \frac{\partial}{\partial s} \frac{e^{-\frac{x^2}{2s^{2H_1}}}}{\sqrt{2\pi s^{2H_1}}} \frac{\partial}{\partial s} \frac{e^{-\frac{s^2}{2t^{2H_2}}}}{\sqrt{2\pi t^{2H_2}}}ds =\left[ \textrm{ by } (\ref{eq:problemBH}) \right]
\end{equation*}
\begin{equation*}
\textcolor{white}{\frac{\partial}{\partial t}p(x,t)}=-2H_2 t^{2H_2 -1}  \int_{0}^{\infty} H_1 s^{2H_1 -1} \frac{\partial^2}{\partial x^2} \frac{e^{-\frac{x^2}{2s^{2H_1}}}}{\sqrt{2\pi s^{2H_1}}} \left( -\frac{s}{t^{2H_2}} \right) \frac{e^{-\frac{s^2}{2t^{2H_2}}}}{\sqrt{2\pi t^{2H_2}}}ds
\end{equation*}
\begin{equation*}
\textcolor{white}{\frac{\partial}{\partial t}p(x,t)}=2H_1 H_2 t^{-1}  \int_{0}^{\infty} s^{2H_1} \frac{\partial^2}{\partial x^2} \frac{e^{-\frac{x^2}{2s^{2H_1}}}}{\sqrt{2\pi s^{2H_1}}} \frac{e^{-\frac{s^2}{2t^{2H_2}}}}{\sqrt{2\pi t^{2H_2}}}ds
\end{equation*}
\begin{equation*}
\textcolor{white}{\frac{\partial}{\partial t}p(x,t)}=-2H_1 H_2 t^{-1} \frac{\partial}{\partial x} \left( x \int_{0}^{\infty} \frac{e^{-\frac{x^2}{2s^{2H_1}}}}{\sqrt{2\pi s^{2H_1}}} \frac{e^{-\frac{s^2}{2t^{2H_2}}}}{\sqrt{2\pi t^{2H_2}}}ds \right)= -H_1 H_2 t^{-1} \frac{\partial}{\partial x} \left( x p \right)
\end{equation*}
and this concludes the proof of theorem.
\begin{flushright}
$\blacksquare$
\end{flushright}

We study in the next theorem the ensemble of solutions to (\ref{eq:firstOrd})
\begin{te}
The first order partial differential equation
\begin{equation}
\frac{t}{H_1 H_2} \frac{\partial u}{\partial t} + x \frac{\partial u}{\partial x}= - u
\label{eq:firstOrd2}
\end{equation}
has a general solution of the form
\begin{equation}
u(x,t)=\frac{1}{x} f\left( \frac{x}{t^{H_1 H_2}} \right), \qquad x \in \mathbb{R}- \{0 \}, \; t>0
\label{solutionFirstOrd}
\end{equation}
with $H_1$, $H_2 \in (0,1)$ and $f \in C^{1}(\mathbb{R})$.
\end{te}
\paragraph{Proof:} The auxiliary equations pertaining to (\ref{eq:firstOrd2}) are
\begin{equation}
H_1 H_2 \frac{dt}{t}=\frac{dx}{x}=-\frac{du}{u}.
\label{eq:lagrang}
\end{equation}
The first couple of members of (\ref{eq:lagrang}) 
\begin{equation*}
H_1 H_2 \frac{dt}{t}=\frac{dx}{x}
\end{equation*}
has solution
\begin{equation*}
c_1 t^{H_1 H_2} = x
\end{equation*}
while equating the second and third term of (\ref{eq:lagrang}) we get
\begin{equation*}
\log x = -\log u +constant
\end{equation*}
and thus
\begin{equation*}
u=\frac{1}{x} c_2= \frac{1}{x} f(c_1) = \frac{1}{x} f\left( \frac{x}{t^{H_1 H_2}} \right).
\end{equation*}
\begin{flushright}
$\blacksquare$
\end{flushright}

\begin{os}
\label{remarkPQ}
\normalfont
We note that the functions of the form 
\begin{equation*}
u(x,t)=\frac{1}{x} f \left( \frac{x}{t^{H_1 H_2}} \right)=\frac{1}{t^{H_1 H_2}} g \left( \frac{x}{t^{H_1 H_2}} \right)
\end{equation*}
with $f,g \in C^{1}(\mathbb{R})$ are solutions to (\ref{eq:firstOrd}) or (\ref{eq:firstOrd2}).\\
The law of the iterated fractional Brownian motion belongs to the class of functions (\ref{solutionFirstOrd}) because it can be written as
\begin{equation*}
\frac{1}{x} \int_{0}^{\infty} \frac{x}{t^{H_1 H_2}} \frac{e^{-\frac{x^2}{2t^{2 H_1 H_2} y^{2H_1}}}}{\sqrt{2\pi y^{2H_1}}} \frac{e^{-\frac{y^2}{2}}}{\sqrt{2\pi}} dy=\frac{1}{x}f\left( \frac{x}{t^{H_1 H_2}} \right)
\end{equation*}
or
\begin{equation*}
\frac{1}{t^{H_1H_2}} \int_{0}^{\infty} \frac{e^{-\frac{x^2}{2t^{2 H_1 H_2}y^{2H_1} }}}{\sqrt{2\pi y^{2H_1}}} \frac{e^{-\frac{y^2}{2}}}{\sqrt{2\pi}} dy =\frac{1}{t^{H_1H_2}}g\left( \frac{x}{t^{H_1 H_2}} \right)
\end{equation*}
Furthermore if we choose $f(z)=\frac{ze^{-\frac{z^2}{2}}}{\sqrt{2\pi}}$ or $g(z)=\frac{e^{-\frac{z^2}{2}}}{\sqrt{2\pi}}$ we see that the Gaussian kernel of Brownian motion belongs to the class of solutions to (\ref{eq:firstOrd}) with $H_1H_2=\frac{1}{2}$. 
\end{os}

\begin{te}
The density of the iterated fractional Brownian motion $I^1_F(t)=B^1_{H_1}(|B^2_{H_2}(t)|)$ for $H_1, H_2 \in (0,1)$ 
\begin{equation}
q(x,t)=2\int_{0}^{\infty} \frac{e^{-\frac{x^2}{2 s^{2H_1}}}}{\sqrt{2 \pi s^{2H_1}}} \frac{e^{-\frac{s^2}{2 t^{2H_2}}}}{\sqrt{2 \pi t^{2H_2}}}ds, \quad x \in \mathbb{R}, \; t>0, \; H_1, H_2 \in (0,1)
\label{densityMMM}
\end{equation}
is a solution to the following second-order p.d.e.
\begin{equation}
(1+H_1 H_2) t \frac{\partial p}{\partial t} + t^2 \frac{\partial^2 p}{\partial t^2} =H_1^2 H_2^2 \left\lbrace 2x\frac{\partial p}{\partial x} + x^2 \frac{\partial^2 p}{\partial x^2} \right\rbrace, \quad  x \in \mathbb{R}, \, t>0, \; H_1, H_2 \in (0,1). 
\label{eqdue}
\end{equation}
\end{te}
We give now two different proofs of this result.

\paragraph{First proof:}
We start by taking the time derivative of \eqref{densityMMM} and by taking into account \eqref{eq:problemBH} we have that
\begin{equation} 
\label{criticalDim}
\frac{\partial p}{\partial t}=2 \int_{0}^{\infty} \frac{e^{-\frac{x^2}{2s^{2H_1}}}}{\sqrt{2 \pi s^{2H_1}}} \frac{\partial}{\partial t} \left( \frac{e^{-\frac{s^2}{2t^{2H_2}}}}{\sqrt{2 \pi t^{2H_2}}} \right) ds 
\end{equation}
\begin{equation*}
\textcolor{white}{\frac{\partial p}{\partial t}} = 2H_2 t^{2H_2-1} \int_{0}^{\infty} \frac{e^{-\frac{x^2}{2s^{2H_1}}}}{\sqrt{2 \pi s^{2H_1}}} \frac{\partial^2}{\partial s^2} \left( \frac{e^{-\frac{s^2}{2t^{2H_2}}}}{\sqrt{2 \pi t^{2H_2}}} \right) ds
\end{equation*}
\begin{equation*}
\textcolor{white}{\frac{\partial p}{\partial t}} = -2H_2 t^{2H_2-1} \int_{0}^{\infty} \frac{\partial}{\partial s} \left( \frac{e^{-\frac{x^2}{2s^{2H_1}}}}{\sqrt{2 \pi s^{2H_1}}} \right) \frac{\partial}{\partial s} \left( \frac{e^{-\frac{s^2}{2t^{2H_2}}}}{\sqrt{2 \pi t^{2H_2}}} \right) ds.
\end{equation*}
In the last step an integration by parts has been carried out.\\ 
The critical point emerges on carring out the second integration by parts which yields
\begin{equation}
\frac{\partial p}{\partial t}=2H_2 t^{2H_2-1}  \int_{0}^{\infty}  \frac{\partial^2}{\partial s^2} \left( \frac{e^{-\frac{x^2}{2s^{2H_1}}}}{\sqrt{2 \pi s^{2H_1}}} \right) \frac{e^{-\frac{s^2}{2t^{2H_2}}}}{\sqrt{2 \pi t^{2H_2}}} ds 
\label{criticalPoint}
\end{equation}
\begin{equation*}
\textcolor{white}{\frac{\partial p}{\partial t}} - 2H_2 t^{2H_2-1} \frac{\partial}{\partial s} \left( \frac{e^{-\frac{x^2}{2s^{2H_1}}}}{\sqrt{2 \pi s^{2H_1}}} \right) \frac{e^{-\frac{s^2}{2t^{2H_2}}}}{\sqrt{2 \pi t^{2H_2}}} \Bigg|_{s=0}^{s=\infty}.
\end{equation*}
By applying once again \eqref{eq:problemBH} we have that
\begin{equation}
\label{criticalPoint2}
\frac{\partial p}{\partial t}=2H_2 t^{2H_2-1}  \int_{0}^{\infty}  \frac{\partial^2}{\partial s^2} \left( \frac{e^{-\frac{x^2}{2s^{2H_1}}}}{\sqrt{2 \pi s^{2H_1}}} \right) \frac{e^{-\frac{s^2}{2t^{2H_2}}}}{\sqrt{2 \pi t^{2H_2}}} ds 
\end{equation}
\begin{equation*}
\textcolor{white}{\frac{\partial p}{\partial t}} +2H_1H_2 t^{2H_2-1} s^{2H_1 -1} \frac{\partial^2}{\partial x^2} \left( \frac{e^{-\frac{x^2}{2s^{2H_1}}}}{\sqrt{2 \pi s^{2H_1}}} \right) \frac{1}{\sqrt{2 \pi t^{2H_2}}} \Bigg|_{s=0}.
\end{equation*}
From \eqref{criticalPoint2} it is clear that for $H_1>1/2$ the second term disappears for $s \to 0^+$ while for $0<H_1 \leq 1/2$ we get a Dirac function as in all previous cases examined above.\\
It is convenient to rewrite \eqref{criticalPoint2} in the following manner
\begin{equation}
\frac{\partial p}{\partial t}=2H_2 t^{2H_2-1}  \int_{0}^{\infty}  \frac{\partial^2}{\partial s^2} \left( \frac{e^{-\frac{x^2}{2s^{2H_1}}}}{\sqrt{2 \pi s^{2H_1}}} \right) \frac{e^{-\frac{s^2}{2t^{2H_2}}}}{\sqrt{2 \pi t^{2H_2}}} ds + \frac{2H_1H_2}{\sqrt{2\pi}} t^{H_2-1} \textbf{1}_{(0 ,1/2]}(H_1) \frac{d^2}{d x^2} \delta(x)
\label{derNOW}
\end{equation}
where the first integral can be developed as follows
\begin{equation}
2H_2 t^{2H_2-1} \int_{0}^{\infty}  \frac{\partial}{\partial s} \left\lbrace  H_1 s^{2H_1-1} \frac{\partial^2}{\partial x^2} \left( \frac{e^{-\frac{x^2}{2s^{2H_1}}}}{\sqrt{2 \pi s^{2H_1}}} \right) \right\rbrace  \frac{e^{-\frac{s^2}{2t^{2H_2}}}}{\sqrt{2 \pi t^{2H_2}}} ds
\label{ulterioriCalc}
\end{equation}
\begin{equation*}
= 2H_2 t^{2H_2-1} \int_{0}^{\infty}  \left\lbrace  H_1(2H_1-1)s^{2H_1-2}\frac{\partial^2}{\partial x^2} \left( \frac{e^{-\frac{x^2}{2s^{2H_1}}}}{\sqrt{2 \pi s^{2H_1}}}\right)  \right .
\end{equation*}
\begin{equation*}
\left . +  H^2_1 s^{4H_1-2} \frac{\partial^4}{\partial x^4} \left( \frac{e^{-\frac{x^2}{2s^{2H_1}}}}{\sqrt{2 \pi s^{2H_1}}} \right) \right\rbrace  \frac{e^{-\frac{s^2}{2t^{2H_2}}}}{\sqrt{2 \pi t^{2H_2}}} ds
\end{equation*}
and can be conveniently rewritten as
\begin{equation*}
2H_1 H_2 t^{H_2-1} \int_{0}^{\infty} s^{2H_1-2} \left[(H_1-1) \frac{\partial^2}{\partial x^2} \left( \frac{e^{-\frac{x^2}{2s^{2H_1}}}}{\sqrt{2 \pi s^{2H_1}}}\right) -  H_1\frac{\partial^2}{\partial x^2} \left( x \frac{\partial}{\partial x} \frac{e^{-\frac{x^2}{2s^{2H_1}}}}{\sqrt{2 \pi s^{2H_1}}} \right) \right] \frac{e^{-\frac{s^2}{2t^{2H_2}}}}{\sqrt{2 \pi}} ds.
\end{equation*}
We note at this point that for $H_1=1/2$ the first integral in the last member of \eqref{ulterioriCalc} disappears and immediately we obtain equation \eqref{pde2}. For $H_1 \neq 1/2$ an additional analysis is necessary and a qualitatively different result is obtained.\\
The time derivative \eqref{derNOW} multiplyed by $t^{1-H_2}$ becomes

\begin{equation*}
t^{1-H_2} \frac{\partial p}{\partial t} = 2H_1 H_2 \int_{0}^{\infty} s^{2H_1-2} \left[(H_1-1) \frac{\partial^2}{\partial x^2} \left( \frac{e^{-\frac{x^2}{2s^{2H_1}}}}{\sqrt{2 \pi s^{2H_1}}}\right) \right .
\end{equation*}
\begin{equation*}
\textcolor{white}{t^{1-H_2} \frac{\partial p}{\partial t}} \left . -  H_1\frac{\partial^2}{\partial x^2} \left( x \frac{\partial}{\partial x} \frac{e^{-\frac{x^2}{2s^{2H_1}}}}{\sqrt{2 \pi s^{2H_1}}} \right) \right] \frac{e^{-\frac{s^2}{2t^{2H_2}}}}{\sqrt{2 \pi}} ds +2\frac{H_1H_2}{\sqrt{2\pi}} \textbf{1}_{(0,\frac{1}{2}]} (H_1) \frac{d^2}{d x^2} \delta(x) .
\end{equation*}
\begin{equation*}
\textcolor{white}{t^{1-H_2} \frac{\partial p}{\partial t}} = 2H_1 H_2 \int_{0}^{\infty} s^{-2} \left[(1-H_1) \frac{\partial}{\partial x} \left( x \frac{e^{-\frac{x^2}{2s^{2H_1}}}}{\sqrt{2 \pi s^{2H_1}}}\right) \right .
\end{equation*}
\begin{equation*}
\textcolor{white}{t^{1-H_2} \frac{\partial p}{\partial t}} \left .  + H_1\frac{\partial^2}{\partial x^2} \left( x^2 \frac{e^{-\frac{x^2}{2s^{2H_1}}}}{\sqrt{2 \pi s^{2H_1}}} \right) \right] \frac{e^{-\frac{s^2}{2t^{2H_2}}}}{\sqrt{2 \pi}} ds +2\frac{H_1H_2}{\sqrt{2\pi}} \frac{d^2}{d x^2} \delta(x) \textbf{1}_{(0,\frac{1}{2}]} (H_1).
\end{equation*}
To get rid of the $s^{-2}$ appearing in the above integral we derive both members w.r. to time $t$ and realise that 
\begin{equation*}
\frac{\partial}{\partial t}\left( t^{1-H_2} \frac{\partial p}{\partial t} \right) = 2\frac{H_1H_2^2}{t^{H_2 +1}} \int_{0}^{\infty}  \left[(1-H_1) \frac{\partial}{\partial x} \left( x \frac{e^{-\frac{x^2}{2s^{2H_1}}}}{\sqrt{2 \pi s^{2H_1}}}\right) \right .
\end{equation*}
\begin{equation*}
\textcolor{white}{\frac{\partial}{\partial t}\left( t^{1-H_2} \frac{\partial p}{\partial t} \right)} \left .  + H_1\frac{\partial^2}{\partial x^2} \left( x^2 \frac{e^{-\frac{x^2}{2s^{2H_1}}}}{\sqrt{2 \pi s^{2H_1}}} \right) \right] \frac{e^{-\frac{s^2}{2t^{2H_2}}}}{\sqrt{2 \pi t^{2H_2}}} ds 
\end{equation*}
\begin{equation*}
\textcolor{white}{\frac{\partial}{\partial t}\left( t^{1-H_2} \frac{\partial p}{\partial t} \right) }=\frac{H_1H_2^2}{t^{H_2 +1}} \left\lbrace (1-H_1) \frac{\partial}{\partial x} \left( x p \right) + H_1\frac{\partial^2}{\partial x^2} \left( x^2 p \right)  \right\rbrace
\end{equation*} 
\begin{equation*}
\textcolor{white}{\frac{\partial}{\partial t}\left( t^{1-H_2} \frac{\partial p}{\partial t} \right) } = \frac{H_1H_2^2}{t^{H_2 +1}} \left\lbrace (1+H_1) \frac{\partial}{\partial x} \left( x p \right) +H_1 \frac{\partial}{\partial x}\left( x^2\frac{\partial p}{\partial x} \right) \right\rbrace 
\end{equation*}
\begin{equation*}
\textcolor{white}{\frac{\partial}{\partial t}\left( t^{1-H_2} \frac{\partial p}{\partial t} \right) } = \frac{H_1H_2^2}{t^{H_2 +1}} \left\lbrace (1+H_1) \frac{\partial}{\partial x} \left( x p \right) +H_1 \left[ 2x \frac{\partial p}{\partial x} + x^2 \frac{\partial^2 p}{\partial x^2} \right] \right\rbrace .
\end{equation*}
In light of theorem \eqref{eq:problemBH}, by writing $-H_1H_2 \frac{\partial}{\partial x}(xp)=t\frac{\partial p}{\partial t}$ and after some calculations, equation \eqref{eqdue} emerges. 
\begin{flushright}
$\blacksquare$
\end{flushright}

\paragraph{Second proof:}
By taking the time derivative of \eqref{eq:firstOrd} multiplyed by $t$ we have that
\begin{align*}
t \frac{\partial }{\partial t}\left( t\frac{\partial p}{\partial t} \right) = & t\frac{\partial p}{\partial t} + t^2 \frac{\partial^2 p}{\partial t^2}\\ 
= & - H_1 H_2 \frac{\partial}{\partial x} \left(x \, t\,\frac{\partial p}{\partial t} \right)\\
= & - H_1 H_2 \frac{\partial }{\partial x} \left( x \left[ -H_1 H_2 \frac{\partial}{\partial x} \left( x p \right) \right] \right)
\end{align*}
\begin{equation*}
\textcolor{white}{t \frac{\partial }{\partial t}\left( t\frac{\partial p}{\partial t} \right)} = H_1^2H_2^2 \left\lbrace \frac{\partial }{\partial x} \left(xp\right) + x\frac{\partial^2 }{\partial x^2} \left(xp\right) \right\rbrace 
\end{equation*}
\begin{equation*}
\textcolor{white}{t \frac{\partial }{\partial t}\left( t\frac{\partial p}{\partial t} \right)} = H_1^2H_2^2  \frac{\partial }{\partial x} \left(xp\right) + H_1^2H_2^2 \left[ 2x \frac{\partial}{\partial x} + x^2 \frac{\partial^2 p}{\partial x^2} \right]
\end{equation*}
\begin{equation*}
\textcolor{white}{t \frac{\partial p}{\partial t}\left( t\frac{\partial p}{\partial t} \right)} = -H_1H_2 t \frac{\partial p}{\partial t} + H_1^2H_2^2 \left[ 2x \frac{\partial p}{\partial x} + x^2 \frac{\partial^2 p}{\partial x^2} \right].
\end{equation*}
In the last step we used again \eqref{eq:firstOrd}.
\begin{flushright}
$\blacksquare$
\end{flushright}

\begin{os}
\normalfont
The equation \eqref{eqdue} displays the structure of the telegraph equation with non-constant coefficients emerging in the case of finite-velocity one-dimensional motions in a non-homogeneous medium (see Ratanov \cite{RT97}). An equation similar to \eqref{eqdue} emerges in the study of hyperbolic random motions of finite velocity in hyperbolic spaces (Orsingher and De Gregorio \cite{ODG07}). In view of remark \ref{remarkPQ} and also by direct calculations it can be proved that the distribution densities of fractional Brownian motion, Brownian motion are solutions to the wave-type equation \eqref{eqdue}.
\end{os}

For the $n$-times fractional Brownian motions we can easily evaluate the even-order moments.
\begin{te}
For the $n$-times iterated fractional Brownian motion 
\begin{equation*}
I^n_F(t)=B^1_{H_1}(|B^2_{H_2}(\ldots |B^{n+1}_{H_{n+1}}(t)| \ldots)|)
\end{equation*}
the explicit form of the moments reads
\begin{equation}
\label{momenti}
E\left\lbrace I^n_F(t) \right\rbrace^{2k}  = \frac{2^{n+1} t^{2k\prod_{j=1}^{n+1}H_j}}{2^{k\sum_{r=0}^{n} \prod_{j=0}^{r}H_j}} \prod_{r=0}^{n} \frac{\Gamma(2k\prod_{j=0}^{r}H_j)}{\Gamma(k\prod_{j=0}^{r}H_j)}, \quad k \in \mathbb{N}, \; n \geq 1
\end{equation}
with $H_0=1$ and $H_j \in (0,1)$ for $j=1,2 \ldots , n$.\\ 
\end{te}
\paragraph{Proof:}
We proceed by recurrence by first evaluating
\begin{equation*}
E\left\lbrace I^1_F(t) \right\rbrace^{2k} = 2 \int_{\mathbb{R}} x^{2k} \int_{0}^{\infty} \frac{e^{-\frac{x^2}{2s^{2H_1}} }}{\sqrt{2 \pi s^{2H_1}}} \frac{e^{-\frac{s^2}{2t^{2H_2}}}}{\sqrt{2 \pi t^{2H_2}}} ds \, dx.
\end{equation*}
Since
\begin{equation*}
2^2 \int_{0}^{\infty} \frac{e^{-\frac{x^2}{2s^{2H_1}}}}{\sqrt{2 \pi s^{2H_1}}} x^{2k} dx=\frac{2^{k+1}}{\sqrt{\pi}} s^{2kH_1} \Gamma\left( k + \frac{1}{2} \right)
\end{equation*}
\begin{equation*}
\textcolor{white}{2^2 \int_{0}^{\infty} \frac{e^{-\frac{x^2}{2s^{H_1}}}}{\sqrt{2 \pi s^{2H_1}}} x^{2k} dx}=2^{2-k} \frac{\Gamma(2k)}{\Gamma(k)} s^{2kH_1},
\end{equation*}
we obtain that
\begin{equation}
\int_0^\infty s^{2kH_1} \frac{e^{-\frac{s^2}{2t^{2H_2}}}}{\sqrt{2 \pi t^{2H_2}}} ds= \frac{2^{kH_1 - 1}}{\sqrt{\pi}} t^{2kH_1H_2} \Gamma\left( kH_1 + \frac{1}{2} \right)
\label{JofTP}
\end{equation}
\begin{equation*}
\textcolor{white}{\int_0^\infty s^{2kH_1} \frac{e^{-\frac{s^2}{2t^{2H_2}}}}{\sqrt{2 \pi t^{2H_2}}} ds} = \frac{ t^{2kH_1H_2}}{2^{kH_1}} \frac{\Gamma(2kH_1)}{\Gamma(kH_1)}
\end{equation*}
and thus
\begin{equation*}
E\left\lbrace I^1_F(t) \right\rbrace^{2k} = 2^{2-k(1+H_1)}  t^{2kH_1H_2} \frac{\Gamma(2k) \Gamma(2kH_1)}{\Gamma(k) \Gamma(kH_1)}.
\end{equation*}
If we assume that result \eqref{momenti} holds for $n-1$ in view of \eqref{JofTP} we immediately have the claimed formula.
\begin{flushright}
$\blacksquare$
\end{flushright}

\begin{os}
\normalfont
The variance of the iterated fractional Brownian motion
\begin{equation}
Var\left\lbrace B^1_{H_1}(|B^2_{H_2}(t)|) \right\rbrace = 2^{1-H_1} \frac{\Gamma(2H_1)}{\Gamma(H_1)} t^{2H_1H_2}
\end{equation}
is an increasing function of the time $t$ which grows more rapidly than the variance of $IBM$
\begin{equation}
Var\left\lbrace B^1(|B^2(t)|) \right\rbrace =\sqrt{\frac{t}{2\pi}}
\end{equation}
for all $H_1$,$H_2$ such that
\[ H_1 H_2 >\frac{1}{2^2}. \]
For $H_1H_2 > 1/2$,  $Var\left\lbrace B^1_{H_1}(|B_{H_2}^2(t)|) \right\rbrace $ increases more rapidly than the variance of standard Brownian motion.
\end{os}

\section{The iterated fractional Brownian motion and modified Bessel functions }
A process related to the iterated fractional Brownian motion $I_F(t)$, $t>0$ is here analysed. Unlike the $n$-times iterated Brownian motion 
\begin{equation}
I^n (t)=B^{1} ( | B^{2} ( \ldots | B^{n+1}(t) | \ldots )| ), \quad t>0
\end{equation}
for which the probability distribution is expressed as a $n$-fold integral, for the process
\begin{equation}
J_{F}^n (t)=B_{H}^{1} ( | B_{H}^{2} ( \ldots | B_H^{n+1}(t) |^\frac{1}{H} \ldots )|^\frac{1}{H} ), \quad t>0
\end{equation}
the probability density can be explicitly expressed in terms of modified Bessel functions. In particular, for $n=1,2,3$ we have the following expressions
\begin{align}
& Pr\left\lbrace J_F^{1}(t) \in dx \right\rbrace = dx \frac{1}{\pi t^H} K_{0} \left( \frac{|x|}{t^H} \right) \label{densityJ1} \\
& Pr\left\lbrace J_F^{2}(t) \in dx \right\rbrace = dx \int_{0}^{\infty} \frac{e^{-\frac{x^2}{2 s^2}}}{\sqrt{2 \pi s^2}} \frac{1}{\pi t^H} K_{0} \left( \frac{|s|}{t^H} \right) ds \label{densityJ2} \\
& Pr\left\lbrace J_F^{3}(t) \in dx \right\rbrace = dx \int_{0}^{\infty} \frac{1}{\pi s} K_{0} \left( \frac{|x|}{s} \right) \frac{1}{\pi t^H} K_{0} \left( \frac{|s|}{t^H} \right)ds.
\end{align}
In general, for $n=2m$ we have that
\begin{equation*}
Pr\left\lbrace J^{2m}_F(t) \in dx \right\rbrace =
\end{equation*}
\begin{equation}
dx \int_{0}^{\infty} ds_1 \ldots \int_{0}^{\infty} ds_m \frac{1}{\pi s_1} K_0 \left( \frac{|x|}{s_1} \right) \ldots \frac{1}{\pi s_{m-1}} K_0 \left( \frac{|s_{m-1}|}{s_m} \right) \frac{1}{\pi t^H} K_0 \left( \frac{|s_{m}|}{t^H} \right)
\end{equation}
and for $n=2m+1$
\begin{align}
Pr\left\lbrace J^{2m+1}_F(t) \in dx \right\rbrace = & dx \int_{0}^{\infty} \frac{e^{-\frac{x^2}{2s_1^2}}}{\sqrt{2 \pi s_1^2}} Pr\left\lbrace J^{2m}_F(t) \in ds_1 \right\rbrace\\
= & \int_{0}^{\infty} Pr\left\lbrace J^{2m}_F(s_n) \in dx \right\rbrace \frac{e^{-\frac{s_n^2}{2t^{2H}}}}{\sqrt{2 \pi t^{2H}}}ds_n .
\end{align}
We shall call $J^n_F(t)$ the weighted iterated fractional Brownian motion.\\
The function $K_\nu (x)$ can be presented in several alternative forms as
\begin{equation}
K_{0}(x)=\int_{0}^{\infty} s^{-1} e^{-\frac{x^2}{4s^2}-s^2} ds, \qquad | \arg x | < \frac{\pi}{4}
\label{functionK0}
\end{equation}
(see p. 119 Lebedev \cite{LE}) or, for $\Re\{ \nu \} > - \frac{1}{2}$
\begin{equation}
K_\nu (x)=\left( \frac{\pi}{2x} \right)^{\frac{1}{2}} \frac{e^{-x}}{\Gamma\left( \nu + \frac{1}{2} \right)} \int_{0}^{\infty} e^{-z} z^{\nu - \frac{1}{2}} \left( 1 + \frac{z}{2x} \right)^{\nu - \frac{1}{2}} dz, \quad | \arg x| < \pi 
\end{equation}
(see p. 140 Lebedev  \cite{LE}) or, for $\nu \neq 0, \; \nu=\pm 1, \pm2, \ldots$
\begin{equation}
K_\nu(x)=\frac{\pi}{2} \frac{I_{-\nu}(x) - I_{\nu}(x)}{\sin \nu \pi}, \quad | \arg x | < \pi
\end{equation}
(see p. 108 Lebedev  \cite{LE}) where
\begin{equation}
I_{\nu}(x)=\sum_{k=0}^{\infty} \left( \frac{x}{2} \right)^{2k + \nu} \frac{1}{k! \Gamma(k+\nu+1)}, \quad |x|<\infty, \; | \arg x | < \pi
\end{equation}
is the Bessel Modified function of the first kind (see p. 108 Lebedev \cite{LE}).\\
For the process
\begin{equation}
\label{eq:processIF}
J_{F}^1(t)=B_{H}^{1} ( | B_{H}^{2}(t)|^\frac{1}{H} ), \quad t>0
\end{equation}
the distribution function, in view of formula (\ref{functionK0}), becomes
\begin{equation}
\label{eq:densityIF}
p^{1}(x,t) = 2 \int_{0}^{\infty} \frac{e^{-\frac{x^2}{2s^2}}}{\sqrt{2\pi s^2}} \frac{e^{-\frac{s^2}{2t^{2H}}}}{\sqrt{2\pi t^{2H}}}ds = \frac{1}{\pi t^H} K_{0}\left( \frac{|x|}{t^H}  \right), \qquad x \in \mathbb{R} \setminus \{ 0 \}, \; t>0, H \in (0,1)
\end{equation}
where the change of variable $s=\sqrt{2z}t^H$ must be introduced.\\
In the spirit of the previous sections we give the partial differential equations governing the distributions \eqref{densityJ1} and \eqref{densityJ2}.
We now state our first result for the process $J^1_F(t)$, $t>0$.
\begin{te}
For the process
\begin{equation}
J_{F}^1(t)=B_{H}^{1} ( | B_{H}^{2} (t)|^\frac{1}{H} ), \quad t>0, H \in (0,1)
\end{equation}
obtained by composing two independent fractional Brownian motions $B^1_H$ and $B^2_H$, the distribution (\ref{eq:densityIF}) for $x\neq 0$, solves the p.d.e.
\begin{equation}
\frac{\partial}{\partial t} p^1=-Ht^{2H-1} \left( 2\frac{\partial^2}{\partial x^2} + x \frac{\partial^3}{\partial x^3} \right) p^1 = -Ht^{2H-1} \frac{\partial^2}{\partial x^2} \left( x\frac{\partial}{\partial x}p^1 \right)
\label{eq:theoremIF}
\end{equation}
for all $x\neq 0$, $t>0$ and $H \in (0,1)$.
\end{te}
\paragraph{First proof:}
We derive equation (\ref{eq:theoremIF}) by two different methods. The first one repeatedly uses the relationship of theorem (\ref{te:BH}) while the second one is based on the properties of the modified Bessel functions $K_{\nu}(z)$.\\
We start by calculating the time derivatives of (\ref{eq:densityIF}) as follows:
\begin{align*}
&\frac{\partial}{\partial t} p^1_H(x,t) = \left[ \textrm{ by } (\ref{eq:problemBH}) \right]\\ 
&=2 \int_{0}^{\infty} \frac{e^{-\frac{x^2}{2s^2}}}{\sqrt{2\pi s^2}} \frac{\partial}{\partial t} p_{H}(s,t) ds\\
&= 2H t^{2H-1} \int_{0}^{\infty} \frac{e^{-\frac{x^2}{2s^2}}}{\sqrt{2\pi s^2}} \frac{\partial^2}{\partial s^2} p_{H}(s,t) ds\\
&= 2H t^{2H-1} \left\lbrace \frac{e^{-\frac{x^2}{2s^2}}}{\sqrt{2\pi s^2}} \frac{\partial}{\partial s} p_{H}(s,t) \Bigg|_{s=0}^{s=\infty}- \int_{0}^{\infty} \frac{\partial}{\partial s} \left( \frac{e^{-\frac{x^2}{2s^2}}}{\sqrt{2\pi s^2}} \right) \frac{\partial}{\partial s} p_{H}(s,t) ds \right\rbrace .
\end{align*}
Since
\begin{equation}
\frac{\partial}{\partial s}p_{H}(x,s)=-\frac{s}{t^{2H}}p_{H}(s,t)
\end{equation}
it is easy to see that both functions $\frac{e^{-\frac{x^2}{2s^2}}}{\sqrt{2\pi s^2}}$ and $\frac{\partial}{\partial s} p_{H}(s,t)$ tend to zero as $s \to \infty$. For $s \to 0$ clearly $\frac{\partial}{\partial s}p_{H}(s,t) \to 0$ while the other term converges to zero because $e^{-\frac{x^2}{2s^2}} \to 0$ as $s \to 0$ for $x \neq 0$.\\
This means that
\begin{equation*}
\frac{\partial}{\partial t}p^1(x,t)=-2H t^{2H-1} \int_{0}^{\infty} \frac{\partial}{\partial s} \left( \frac{e^{-\frac{x^2}{2s^2}}}{\sqrt{2\pi s^2}} \right) \frac{\partial}{\partial s}p_{H}(s,t)ds
\end{equation*}
\begin{equation*}
\textcolor{white}{\frac{\partial}{\partial t}p^1(x,t)} = - 2Ht^{2H-1}\left\lbrace \frac{\partial}{\partial s} \left( \frac{e^{-\frac{x^2}{2s^2}}}{\sqrt{2\pi s^2}} \right)p_{H}(s,t)\Bigg|_{s=0}^{s=\infty}- \int_{0}^{\infty}  \frac{\partial^2}{\partial s^2} \left( \frac{e^{-\frac{x^2}{2s^2}}}{\sqrt{2\pi s^2}} \right) p_{H}(s,t)ds \right\rbrace
\end{equation*}
\begin{equation*}
\textcolor{white}{\frac{\partial}{\partial t}p^1(x,t)}= \left[ \textrm{ by } (\ref{eq:problemB1}) \right]
\end{equation*}
\begin{equation*}
\textcolor{white}{\frac{\partial}{\partial t}p^1(x,t)}=2Ht^{2H -1}\left\lbrace \int_{0}^{\infty} \left[ \frac{\partial^2}{\partial x^2} \left( \frac{e^{-\frac{x^2}{2s^2}}}{\sqrt{2\pi s^2}} \right) + s^2 \frac{\partial^4}{\partial x^4} \left( \frac{e^{-\frac{x^2}{2s^2}}}{\sqrt{2\pi s^2}} \right) \right] p_{H}(s,t) ds \right\rbrace .
\end{equation*}
\begin{equation*}
\textcolor{white}{\frac{\partial}{\partial t}p^1(x,t)}= Ht^{2H-1} \left\lbrace 2\int_{0}^{\infty} \left[ \frac{\partial^2}{\partial x^2} \left( \frac{e^{-\frac{x^2}{2s^2}}}{\sqrt{2\pi s^2}} \right) + s^2 \frac{\partial^3}{\partial x^3} \left( -\frac{x}{s^2}  \frac{e^{-\frac{x^2}{2s^2}}}{\sqrt{2\pi s^2}} \right) \right]  p_{H}(s,t) ds \right\rbrace 
\end{equation*}
\begin{equation*}
\textcolor{white}{\frac{\partial}{\partial t}p^1(x,t)}= Ht^{2H-1} \left\lbrace \frac{\partial^2}{\partial x^2}p^1(x,t) - \frac{\partial^3}{\partial x^3} \left[ x p^1(x,t) \right] \right\rbrace
\end{equation*}
\begin{equation*}
\textcolor{white}{\frac{\partial}{\partial t}p^1(x,t)}=- Ht^{2H-1} \left\lbrace 2\frac{\partial^2}{\partial x^2} + x\frac{\partial^3}{\partial x^3}  \right\rbrace p^1(x,t) .
\end{equation*}
In the previous steps it is important to consider that
\begin{equation*}
\frac{\partial}{\partial s} \left( \frac{e^{-\frac{x^2}{2s^2}}}{\sqrt{2\pi s^2}} \right) = \frac{\partial}{\partial s} \left( \frac{1}{2\pi} \int_{-\infty}^{\infty} e^{-i\alpha x - \frac{1}{2}\alpha^2 s^2} d\alpha \right) = \frac{1}{2\pi} \int_{-\infty}^{\infty} e^{-i\alpha x - \frac{1}{2}\alpha^2 s^2} (-\alpha^2 s) d\alpha
\end{equation*}
where the inversion of derivative and integral is justified since 
\[ \Big| e^{-i\alpha x - \frac{1}{2}\alpha^2 s^2} \Big| = e^{-\frac{1}{2}\alpha^2 s^2} \]
possesses finite integral. For the same reason the limit w.r. to $s$ can be brought inside the integral and
\[ \lim_{s \to 0^+} \frac{\partial}{\partial s} \left( \frac{e^{-\frac{x^2}{2s^2}}}{\sqrt{2\pi s^2}} \right) p_{H}(s,t) = 0 . \] 
\begin{flushright}
$\blacksquare$
\end{flushright}

\paragraph{Second proof:}
This proof is based on the following properties of the modified Bessel function (see p. 110 Lebedev \cite{LE})
\begin{equation}
\frac{d}{dz}K_{\nu}(z)=\frac{\nu}{z}K_{\nu}(z) - K_{\nu + 1}(z) 
\label{eq:BesselK1}
\end{equation}
and
\begin{equation}
K_{\nu +1}(z) = K_{\nu - 1}(z) + 2\frac{\nu}{z}K_{\nu}(z).
\label{eq:BesselK2}
\end{equation}
We perform our calculations on the representation (\ref{eq:densityIF}) of the distribution of the iterated fractional Brownian motion. We need to evaluate the partial derivatives appearing in equation (\ref{eq:theoremIF}) for $x \in \mathbb{R}^{+}$, for the derivative w.r. to time we have that
\begin{equation}
\label{eq:timederivativeK0}
\frac{\partial}{\partial t}\left[ \frac{1}{\pi t^H} K_{0}\left( \frac{x}{t^H} \right) \right] = - \frac{H}{\pi t^{H+1}} K_{0}\left( \frac{x}{t^H} \right) - \frac{Hx}{\pi t^{2H+1}}\frac{d}{dz}K_{0}(z)\Bigg|_{z=\frac{x}{t^H}}
\end{equation}
\begin{equation*}
\textcolor{white}{\frac{\partial}{\partial t}\left[ \frac{1}{\pi t^H} K_{0}\left( \frac{x}{t^H} \right) \right] }= \left[ \textrm{ by } (\ref{eq:BesselK1}) \right] 
\end{equation*}
\begin{equation*}
\textcolor{white}{\frac{\partial}{\partial t}\left[ \frac{1}{\pi t^H} K_{0}\left( \frac{x}{t^H} \right) \right] }= \frac{Hx}{\pi t^{2H+1}} K_{1}\left( \frac{x}{t^{H}} \right) -\frac{H}{\pi t^{H+1}} K_{0}\left( \frac{x}{t^H} \right).
\end{equation*}
We now pass to the partial derivatives with respect to space. We have successively that
\begin{equation}
\frac{\partial}{\partial x}\left[ \frac{1}{\pi t^H} K_{0}\left( \frac{x}{t^H} \right) \right]= -\frac{1}{\pi t^{2H}}K_{1}\left( \frac{x}{t^H} \right),
\end{equation}
\begin{equation}
\frac{\partial^2}{\partial x^2} \left[ \frac{1}{\pi t^H} K_{0}\left( \frac{x}{t^H} \right) \right]=\left[ \textrm{ by } (\ref{eq:BesselK1}) \right]
\end{equation}
\begin{equation*}
\textcolor{white}{\frac{\partial^2}{\partial x^2} \left[ \frac{1}{\pi t^H} K_{0}\left( \frac{x}{t^H} \right) \right]} =-\frac{1}{\pi t^{3H}} \frac{d}{dz}K_{1}(z) \Bigg|_{z=\frac{x}{t^H}}
\end{equation*}
\begin{equation*}
\textcolor{white}{\frac{\partial^2}{\partial x^2} \left[ \frac{1}{\pi t^H} K_{0}\left( \frac{x}{t^H} \right) \right]}=-\frac{1}{\pi t^{3H}} \left\lbrace \frac{t^H}{x} K_{1} \left( \frac{x}{t^H} \right) - K_{2}\left(\frac{x}{t^H} \right) \right\rbrace
\end{equation*}
\begin{equation*}
\textcolor{white}{\frac{\partial^2}{\partial x^2} \left[ \frac{1}{\pi t^H} K_{0}\left( \frac{x}{t^H} \right) \right]} =-\frac{1}{\pi t^{3H}} \left\lbrace -K_{0}\left( \frac{x}{t^H} \right) - \frac{t^H}{x}K_{1}\left( \frac{x}{t^H} \right) \right\rbrace .
\end{equation*}
Some more calculations are necessary for the third order partial derivative which reads
\begin{equation}
\frac{\partial^3}{\partial x^3}\left[ \frac{1}{\pi t^H} K_{0}\left( \frac{x}{t^H} \right) \right]=-\frac{1}{\pi t^{3H}} \left\lbrace \frac{1}{t^H}K_{1}\left(\frac{x}{t^H} \right) + \frac{t^H}{x^2} K_{1} \left( \frac{x}{t^H} \right) \right . 
\end{equation}
\begin{equation*}
\textcolor{white}{\frac{\partial^3}{\partial x^3}\left[ \frac{1}{\pi t^H} K_{0}\left( \frac{x}{t^H} \right) \right]=} \left . - \frac{t^H}{x} \left[ \frac{t^H}{x} K_{1}\left( \frac{x}{t^H} \right) - K_{2}\left( \frac{x}{t^H} \right) \right]  \frac{1}{t^H} \right\rbrace
\end{equation*}
\begin{equation*}
\textcolor{white}{\frac{\partial^3}{\partial x^3}\left[ \frac{1}{\pi t^H} K_{0}\left( \frac{x}{t^H} \right) \right] }=-\frac{1}{\pi t^{3H}} \left\lbrace \frac{1}{t^H}K_{1}\left(\frac{x}{t^H} \right) + \frac{1}{x} K_{2}\left( \frac{x}{t^H} \right) \right\rbrace
\end{equation*}
\begin{equation*}
\textcolor{white}{\frac{\partial^3}{\partial x^3}\left[ \frac{1}{\pi t^H} K_{0}\left( \frac{x}{t^H} \right) \right]}=\left[ \textrm{ by } (\ref{eq:BesselK2}) \right]
\end{equation*}
\begin{equation*}
\textcolor{white}{\frac{\partial^3}{\partial x^3}\left[ \frac{1}{\pi t^H} K_{0}\left( \frac{x}{t^H} \right) \right]}=-\frac{1}{\pi t^{3H}} \left\lbrace \frac{1}{x} K_{0} \left( \frac{x}{t^H} \right) + \left( 2 \frac{t^H}{x^2} + \frac{1}{t^H} \right) K_{1} \left( \frac{x}{t^H} \right) \right\rbrace. 
\end{equation*}
Therefore, by combining the previous results one obtains
\begin{equation*}
\left( 2\frac{\partial^2}{\partial x^2} + x \frac{\partial^3}{\partial x^3} \right) \left[ \frac{1}{\pi t^H} K_{0}\left( \frac{x}{t^H} \right) \right]=-\frac{1}{\pi t^{3H}} \left\lbrace - K_{0}\left( \frac{x}{t^H}\right) + \frac{x}{t^H}K_{1} \left( \frac{x}{t^H}\right) \right\rbrace
\end{equation*}
and thus
\begin{equation*}
-H t^{2H-1} \left( 2\frac{\partial^2}{\partial x^2} + x \frac{\partial^3}{\partial x^3} \right) \left[ \frac{1}{\pi t^H} K_{0}\left( \frac{x}{t^H} \right) \right]=\frac{H}{\pi t^{H+1}} \left\lbrace - K_{0}\left( \frac{x}{t^H}\right) + \frac{x}{t^H}K_{1} \left( \frac{x}{t^H}\right) \right\rbrace
\end{equation*}
which coincides with the time derivative (\ref{eq:timederivativeK0}).\\
This concludes the second proof.
\begin{flushright}
$\blacksquare$
\end{flushright}

We now examine the twice iterated fractional Brownian motion and show that its probability distribution solves a fourth-order p.d.e. 
\begin{te}
The process $ J^{2}_{F}(t)=B_{H}^{1} ( | B_{H}^{2}( | B_{H}^{3}(t) |^{\frac{1}{H}} )|^{\frac{1}{H}} )$
has a probability density
\begin{equation}
\label{eq:densityI2F}
p^{2}(x,t)=4\int_{0}^{\infty} \frac{e^{-\frac{x^2}{2s^{2}}}}{\sqrt{2 \pi s^{2}}} \int_{0}^{\infty} \frac{e^{-\frac{s^2}{2z^{2}}}}{\sqrt{2 \pi z^{2}}} \frac{e^{-\frac{z^2}{2t^{2H}}}}{\sqrt{2 \pi t^{2H}}} dzds
\end{equation}
\begin{equation*}
\textcolor{white}{p^{2}(x,t;H)} = 2\int_{0}^{\infty} \frac{e^{-\frac{x^2}{2s^{2}}}}{\sqrt{2 \pi s^{2}}} \frac{1}{\pi t^{H}} K_{0}\left( \frac{s}{t^H} \right) ds
\end{equation*}
\begin{equation*}
\textcolor{white}{p^{2}(x,t;H)} =2\int_{0}^{\infty} \frac{e^{-\frac{x^2}{2s^{2}}}}{\sqrt{2 \pi s^{2}}} p^{1}(s,t)ds, \quad x \in \mathbb{R}, \; t>0, \; H \in (0,1)
\end{equation*}
satisfying the following fourth-order p.d.e.
\begin{equation}
\frac{\partial}{\partial t} p^{2}(x,t)=Ht^{2H-1}\left\lbrace 4\frac{\partial^2}{\partial x^2} + 5x\frac{\partial^3}{\partial x^3} + x^2 \frac{\partial^4}{\partial x^4} \right\rbrace p^{2}(x,t) 
\end{equation}
for $x\neq 0$, $t>0$.
\end{te}
\paragraph{Proof:}
We start our proof by writing down the time derivative w.r. to time $t$ of (\ref{eq:densityI2F}) and making use of the result \eqref{eq:theoremIF} we get that
\begin{equation}
\label{eq:I2Ftimederivative}
\frac{\partial}{\partial t}p^{2}(x,t)=2\int_{0}^{\infty} \frac{e^{-\frac{x^2}{2s^{2}}}}{\sqrt{2 \pi s^{2}}} \frac{\partial}{\partial t} p^{1}(s,t)
\end{equation}
\begin{equation*}
\textcolor{white}{\frac{\partial}{\partial t}p^{2}(x,t)} =-2Ht^{2H-1} \int_{0}^{\infty} \frac{e^{-\frac{x^2}{2s^{2}}}}{\sqrt{2 \pi s^{2}}} \left[ 2\frac{\partial^2}{\partial s^2} + s\, \frac{\partial^3}{\partial s^3} \right] p^{1}(s,t) ds .
\end{equation*}
By applying a couple of integrations by parts we have that
\begin{equation}
\label{eq:I2Fparts1}
\int_{0}^{\infty} \frac{e^{-\frac{x^2}{2s^{2}}}}{\sqrt{2 \pi s^{2}}} \frac{\partial^2}{\partial s^2}p^{1}(s,t) ds=\int_{0}^{\infty} \frac{\partial^2}{\partial s^2} \left( \frac{e^{-\frac{x^2}{2s^{2}}}}{\sqrt{2 \pi s^{2}}} \right) p^{1}(s,t) ds
\end{equation}
while by applying three integrations by parts we obtain
\begin{equation}
\int_{0}^{\infty} \frac{e^{-\frac{x^2}{2s^{2}}}}{\sqrt{2 \pi s^{2}}} \, s\, \frac{\partial^3}{\partial s^3}p^{1}(s,t) ds =- \int_{0}^{\infty} \frac{\partial^3}{\partial s^3} \left( s\, \frac{e^{-\frac{x^2}{2s^{2}}}}{\sqrt{2 \pi s^{2}}} \right) p^{1}(s,t) ds
\label{eq:I2Fparts2}
\end{equation}
\begin{equation*}
\textcolor{white}{\int_{0}^{\infty} \frac{e^{-\frac{x^2}{2s^{2}}}}{\sqrt{2 \pi s^{2}}} \, s\, \frac{\partial^3}{\partial s^3}p^{1}(s,t) ds} =-\int_{0}^{\infty} \frac{\partial^2}{\partial s^2} \left\lbrace \frac{e^{-\frac{x^2}{2s^{2}}}}{\sqrt{2 \pi s^{2}}} + s\frac{\partial}{\partial s} \left( \frac{e^{-\frac{x^2}{2s^{2}}}}{\sqrt{2 \pi s^{2}}} \right) \right\rbrace p^{1}(s,t)ds 
\end{equation*}
In light of (\ref{eq:I2Fparts1}) and (\ref{eq:I2Fparts2}) the time derivative (\ref{eq:I2Ftimederivative}) becomes 
\begin{equation*}
\frac{\partial}{\partial t}p^{2}(x,t)=
\end{equation*}
\begin{equation*}
-2Ht^{2H-1}\int_{0}^{\infty} \left\lbrace  \frac{\partial^2}{\partial s^2} \left( \frac{e^{-\frac{x^2}{2s^{2}}}}{\sqrt{2 \pi s^{2}}} \right) - \frac{\partial^2}{\partial s^2} \left[ s\, \frac{\partial}{\partial s} \left( \frac{e^{-\frac{x^2}{2s^{2}}}}{\sqrt{2 \pi s^{2}}} \right) \right] \right\rbrace  p^{1}(s,t) ds=
\end{equation*}
\begin{equation*}
-2Ht^{2H-1} \int_{0}^{\infty} \left\lbrace  \frac{\partial^2}{\partial s^2} \left( \frac{e^{-\frac{x^2}{2s^{2}}}}{\sqrt{2 \pi s^{2}}} \right) - \frac{\partial}{\partial s}\left[ \frac{\partial}{\partial s} \left( \frac{e^{-\frac{x^2}{2s^{2}}}}{\sqrt{2 \pi s^{2}}} \right) + s \frac{\partial^2}{\partial s^2} \left( \frac{e^{-\frac{x^2}{2s^{2}}}}{\sqrt{2 \pi s^{2}}} \right) \right] \right\rbrace  p^{1}(s,t) ds=
\end{equation*}
\begin{equation}
2Ht^{2H-1} \int_{0}^{\infty} \left\lbrace  \frac{\partial^2}{\partial s^2} \left( \frac{e^{-\frac{x^2}{2s^{2}}}}{\sqrt{2 \pi s^{2}}} \right) + s\, \frac{\partial^3}{\partial s^3} \left( \frac{e^{-\frac{x^2}{2s^{2}}}}{\sqrt{2 \pi s^{2}}} \right) \right\rbrace  p^{1}(s,t) ds.
\end{equation}
For simplicity, we now write
\[ q=q(x,s)=\frac{e^{-\frac{x^2}{2s^{2}}}}{\sqrt{2 \pi s^{2}}}. \]
By repeated use of \eqref{eq:problemB1} we have that 
\begin{equation*}
\frac{\partial^2}{\partial s^2}q= \frac{\partial}{\partial s}\left[ s \frac{\partial^2}{\partial x^2}q \right] = \frac{\partial^2}{\partial x^2}q + s^2 \frac{\partial^4 q}{\partial x^4}
\end{equation*}
and
\begin{equation*}
\frac{\partial^3 q}{\partial s^3}=\frac{\partial}{\partial s}\left[ \frac{\partial^2 q}{\partial x^2} + s^2 \frac{\partial^4 q}{\partial x^4} \right]=3s\frac{\partial^4 q}{\partial x^4} + s^3\frac{\partial^6 q}{\partial x^6}.
\end{equation*}
Hence, by combining the previous results we can rewrite the time derivative \eqref{eq:I2Ftimederivative} as follows
\begin{equation}
\frac{\partial}{\partial t}p^{2}(x,t)=2Ht^{2H-1} \int_{0}^{\infty} \left\lbrace  \frac{\partial^2}{\partial x^2}q(x,s) + 4s^2 \frac{\partial^4}{\partial x^4}q(x,s) + s^4 \frac{\partial^6}{\partial x^6}q(x,s)  \right\rbrace  p^{1}(s,t) ds.
\end{equation}
In order to handle the operator
\begin{equation}
\label{oper}
\left\lbrace  \frac{\partial^2}{\partial x^2} + 4s^2 \frac{\partial^4}{\partial x^4} + s^4 \frac{\partial^6}{\partial x^6}  \right\rbrace q(x,s)
\end{equation}
it will be useful to note that
\begin{equation*}
\frac{\partial}{\partial x} q(x,s) = -\frac{x}{s^2}q(x,s)
\end{equation*}
and
\begin{equation*}
\frac{\partial^2}{\partial x^2} q(x,s) = \left( -\frac{1}{s^2} + \frac{x^2}{s^4} \right)q(x,s),
\end{equation*}
this implies that the operator \eqref{oper} becomes
\begin{equation}
\label{eq:deraux1}
\frac{\partial^2}{\partial x^2}q + 4s^2 \frac{\partial^4}{\partial x^4} q + s^4 \frac{\partial^6}{\partial x^6}q
\end{equation}
\begin{equation*}
=\frac{\partial^2}{\partial x^2}q + 4s^2 \frac{\partial^3}{\partial x^3}\left(- \frac{x}{s^2} \right)q + s^4 4 \frac{\partial^4}{\partial x^4}\left( -\frac{1}{s^2} + \frac{x^2}{s^4} \right)q
\end{equation*}
\begin{equation*}
=\frac{\partial^2}{\partial x^2}q - 4\frac{\partial^3}{\partial x^3}\left( x\, q \right) - s^2 \frac{\partial^4}{\partial x^4}q + \frac{\partial^4}{\partial x^4} \left( x^2 \, q \right)
\end{equation*}
\begin{equation*}
=\frac{\partial^2}{\partial x^2}q - 4\frac{\partial^3}{\partial x^3}\left( x\, q \right) + \frac{\partial^3}{\partial x^3}\left( x\, q \right) + \frac{\partial^4}{\partial x^4}\left( x^2\, q \right)
\end{equation*}
\begin{equation*}
=\frac{\partial^2}{\partial x^2}q - 3\frac{\partial^3}{\partial x^3}\left( x\, q \right) + \frac{\partial^4}{\partial x^4}\left( x^2\, q \right).
\end{equation*}
By some calculations we have that
\begin{equation}
\label{eq:deraux2}
\frac{\partial^3}{\partial x^3}\left( x\, q \right)=3 \frac{\partial^2}{\partial x^2}q + x\frac{\partial^3}{\partial x^3} q
\end{equation}
and
\begin{equation}
\label{eq:deraux3}
\frac{\partial^4}{\partial x^4}\left( x^2\, q \right)= \frac{\partial^3}{\partial x^3}\left( 2x\, q  + x^2 \frac{\partial}{\partial x}q \right)
\end{equation}
\begin{equation*}
\textcolor{white}{\frac{\partial^4}{\partial x^4}\left( x^2\, q \right)}=2 \left[ 3\frac{\partial^2}{\partial x^2}q + x\frac{\partial^3}{\partial x^3}q \right] + \frac{\partial^2}{\partial x^2}\left( 2x \frac{\partial}{\partial x}q + x^2 \frac{\partial^2}{\partial x^2}q \right)
\end{equation*}
\begin{equation*}
\textcolor{white}{\frac{\partial^4}{\partial x^4}\left( x^2\, q \right)}=2 \left[ 3\frac{\partial^2}{\partial x^2}q + x\frac{\partial^3}{\partial x^3}q \right] +2\left[ \frac{\partial^2}{\partial x^2}q + \frac{\partial}{\partial x} \left( x\, \frac{\partial^2}{\partial x^2}q \right) \right]
\end{equation*}
\begin{equation*}
\textcolor{white}{\frac{\partial^4}{\partial x^4}\left( x^2\, q \right)}+ \frac{\partial}{\partial x} \left[ 2x\frac{\partial^2}{\partial x^2}q + x^2 \frac{\partial^3}{\partial x^3}q \right]
\end{equation*}
\begin{equation*}
\textcolor{white}{\frac{\partial^4}{\partial x^4}\left( x^2\, q \right)}= 12\frac{\partial^2}{\partial x^2}q + 8x\frac{\partial^3}{\partial x^3}q + x^2 \frac{\partial^4}{\partial x^4}q,
\end{equation*}
by inserting (\ref{eq:deraux2}) and \eqref{eq:deraux3} into (\ref{eq:deraux1}) we have that
\begin{align*}
&\frac{\partial^2}{\partial x^2}q + 4s^2 \frac{\partial^4}{\partial x^4} q + s^4 \frac{\partial^6}{\partial x^6}q=4\frac{\partial^2}{\partial x^2}q + 5x\frac{\partial^3}{\partial x^3}q + x^2 \frac{\partial^4}{\partial x^4}q
\end{align*} 
and this leads to the claimed result.
\begin{flushright}
$\blacksquare$
\end{flushright}

Let's consider now the process
\begin{equation*} 
J^n_F(t)=B_{H_1}^{1} (| B_{H_2}^{2}(| \ldots | B^{n+1}_{H_{n+1}}(t) |^{\frac{1}{H_{n}}} \ldots |^\frac{1}{H_2} ) |^\frac{1}{H_1}), \quad t>0, \; H_j \in (0,1), \; j=1,2,\ldots , n.
\end{equation*}
We have already seen a particular case in \eqref{JnDistribution1} and \eqref{JnDistribution2} where $H_1=H_2=\ldots =H_{n+1}=H$.  We show that the following equivalence in distribution holds
\begin{equation} 
J^n_F(t) =B_{H_1}^{1} (| B_{H_2}^{2}(| \ldots | B^{n+1}_{H_{n+1}}(t) |^{\frac{1}{H_{n}}} \ldots |^\frac{1}{H_2} ) |^\frac{1}{H_1})\stackrel{i.d.}{=} B^{1} (| B^{2}( \ldots | B^{n+1}_{H_{n+1}}(t) |^2 \ldots ) |^2) .
\label{process:JnF}
\end{equation}
It must be also noted that the unique relevant Hurst parameter $H$ in $J^n_F(t)$ is that in $B^{n+1}_{H_{n+1}}(t)$. For this reason and for the sake of simplicity we assume throughout the paper that $H_j=H$ for $j=1,2,\ldots, n+1$.\\
We can state the following result on the factorization of the iterated fractional Brownian motion.
\begin{te}
For the process 
\begin{equation}
J^{n-1}_F (t)= B_H^{1} (| B_H^{2}( \ldots | B^{n}_{H}(t) |^{\frac{1}{H}} \ldots ) |^\frac{1}{H}), \quad  t>0, \; H \in (0,1)  
\end{equation}
with $B^j_H$, $j=1,2,\ldots , n$ independent fractional Brownian motions, the following equalities in distribution hold
\begin{equation}
J^{n-1}_F (t) \stackrel{i.d.}{=} B^{1} (| B^{2}( \ldots | B^{n}_{H}(t) |^2 \ldots ) |^2) \stackrel{i.d.}{=} \prod_{i=1}^{n} B^i_{\frac{H}{n}} (t), \quad t>0, \; H \in (0,1)
\label{proofuno}
\end{equation}
where $B^j(t)$, $t>0$, $j=1,2,\ldots , n$ are independent standard Brownian motions.
\end{te}
\paragraph{Proof:}
The proof of the first equality in (\ref{proofuno}) is a trivial matter.  Indeed the equivalence in distribution follows from the property $B_H(t) \stackrel{i.d.}{=} B(t^{2H})$ for each $t$. For the proof of the second equality in (\ref{proofuno}) we proceed by evaluating the Mellin transform of the density of $|J^{n-1}_F(t)|$, $t>0$ as
\begin{equation*}
\mathcal{M}_{n-1}(\alpha, t)=\int_{0}^{\infty} x^{\alpha -1} 2^{n} \int_{0}^\infty ds_1 \ldots \int_0^\infty ds_{n-1} \frac{e^{-\frac{x^{2}}{2s_1^2}}}{\sqrt{2 \pi s_1^2}} \frac{e^{-\frac{s^{2}_1}{2s_2^2}}}{\sqrt{2 \pi s_2^2}} \ldots \frac{e^{-\frac{s^{2}_{n-1}}{2t^{2H}}}}{\sqrt{2 \pi t^{2H}}} dx
\end{equation*}
\begin{equation*}
\textcolor{white}{\mathcal{M}_{n-1}(\alpha, t)} = \left[ \frac{2^{\frac{\alpha}{2}} \Gamma\left( \frac{\alpha}{2} \right)}{\sqrt{2\pi}} \right]^n t^{H(\alpha -1)},\quad  \Re \{ \alpha \} >0.
\end{equation*}
The Fourier transform of the distribution of $J^{n-1}_F(t)$ can be written as
\begin{equation}
\label{FourierProd}
\mathcal{F}_{n-1}(\beta, t)= \int_{\mathbb{R}} e^{i \beta x} p_{n-1}(x,t)dx
\end{equation}
\begin{equation*}
\textcolor{white}{\hat{p}_{n-1}(\beta, t)} = \sum_{k =0}^{\infty} \frac{(i\beta)^k}{k!} \int_{\mathbb{R}} x^k p_{n-1}(x,t)dx
\end{equation*}
\begin{equation*}
\textcolor{white}{\hat{p}_{n-1}(\beta, t)} = \sum_{k =0}^{\infty} \frac{(i\beta)^{2k}}{(2k)!} \int_{\mathbb{R}} x^{2k} p_{n-1}(x,t)dx
\end{equation*}
\begin{equation*}
\textcolor{white}{\hat{p}_{n-1}(\beta, t)} = \sum_{k =0}^{\infty} \frac{(i\beta)^{2k}}{(2k)!} \left[ \frac{2^{k+\frac{1}{2}} \Gamma\left( k + \frac{1}{2} \right)}{\sqrt{2\pi}} \right]^n t^{2Hk}.
\end{equation*}
We now evaluate the characteristic function of 
\begin{equation*}
\prod_{j=1}^{n} B^j_{\frac{H}{n}} (t)
\end{equation*} 
which becomes
\begin{equation*}
E e^{i \beta \prod_{j=1}^{n} B^j_{\frac{H}{n}} (t)} =\sum_{k =0}^{\infty} \frac{(i\beta)^{k}}{k!} E\left\lbrace \prod_{j=1}^{n} B^j_{\frac{H}{n}} (t) \right\rbrace 
\end{equation*}
\begin{equation*}
\textcolor{white}{E e^{i \beta \prod_{j=1}^{n} B^j_{\frac{H}{n}} (t)}}= \sum_{k =0}^{\infty} \frac{(i\beta)^{k}}{k!}  \prod_{j=1}^{n} E\left\lbrace B^j_{\frac{H}{n}} (t) \right\rbrace^k 
\end{equation*}
\begin{equation*}
\textcolor{white}{E e^{i \beta \prod_{j=1}^{n} B^j_{\frac{H}{n}} (t)}}= \sum_{k =0}^{\infty} \frac{(i\beta)^{2k}}{(2k)!}  \prod_{j=1}^{n} E\left\lbrace B^j_{\frac{H}{n}} (t) \right\rbrace^{2k} 
\end{equation*}
\begin{equation*}
\textcolor{white}{E e^{i \beta \prod_{j=1}^{n} B^j_{\frac{H}{n}} (t)}}=\sum_{k =0}^{\infty} \frac{(i\beta)^{2k}}{(2k)!} \prod_{j=1}^{n} \left[ \frac{2^{k+\frac{1}{2}} \Gamma\left(k+ \frac{1}{2} \right)}{\sqrt{2\pi}} \right] t^{\frac{H}{n} 2k}
\end{equation*}
\begin{equation*}
\textcolor{white}{E e^{i \beta \prod_{j=1}^{n} B^j_{\frac{H}{n}} (t)}}=\sum_{k =0}^{\infty} \frac{(i\beta)^{2k}}{(2k)!} \left[ \frac{2^{k+\frac{1}{2}} \Gamma\left(k+ \frac{1}{2} \right)}{\sqrt{2\pi}} \right]^{n} t^{2Hk}
\end{equation*}
and this coincides with \eqref{FourierProd}. This confirms the equivalences in distribution of \eqref{proofuno}.
\begin{flushright}
$\blacksquare$
\end{flushright}
\begin{os}
\normalfont
The result of the previous theorem permits us to write down the following special factorization
\begin{equation*}
B^1_H(|B^2_H(t)|^\frac{1}{H}) \stackrel{i.d.}{=} B^1_{\frac{H}{2}}(t) B^2_{\frac{H}{2}}(t) \stackrel{i.d.}{=} B^1(|B_H^2(t)|^2), \quad t>0, \; H \in (0,1).
\end{equation*}
\end{os}

\section{The iterated Cauchy process}
The Cauchy process has been investigated from many viewpoints and some authors have examined the structure of the jumps, some others have concentrated their analysis on the behaviour of sample paths and others have obtained a number of its distributional properties. We here consider the composition of independent Cauchy processes and its connection with wave equations.\\
The Cauchy process $C(t)$, $t >0$ has a law 
\[ p(x,t)=\frac{t}{\pi (t^2 + x^2)} \quad x \in \mathbb{R}, \; t>0 \]
which satisfies the Laplace equation
\begin{equation}
\frac{\partial^2 p}{\partial t^2} + \frac{\partial^2 p}{\partial x^2}=0 \quad x \in \mathbb{R}, \; t>0
\label{Cuno}
\end{equation} 
and also the space-fractional equation
\begin{equation}
\frac{\partial p}{\partial t} = -\frac{\partial p}{\partial |x|} \quad x \in \mathbb{R}, \; t>0.
\label{Cdue}
\end{equation}
While (\ref{Cuno}) can be checked straightforwardly some additional care is needed for (\ref{Cdue}).\\
The first order spatial derivative appearing in (\ref{Cdue}) is defined as
\begin{equation}
\left( \frac{\partial f}{\partial |x|} \right)(x)= \frac{1}{\pi} \frac{d}{dx} \left\lbrace \int_{-\infty}^{x} \frac{f(s)}{(x-s)} ds - \int_{x}^{\infty} \frac{f(s)}{(s-x)}ds \right\rbrace .
\label{OperatoreR}
\end{equation}
We now evaluate the Fourier transform of (\ref{OperatoreR}).
\begin{equation*}
\int_{\mathbb{R}} e^{i\beta x} \left( \frac{\partial f}{\partial |x|} \right)(x) dx =\frac{1}{\pi}  \int_{\mathbb{R}} e^{i\beta x} \frac{d}{dx} \left\lbrace \int_{-\infty}^{x} \frac{f(s)}{(x-s)} ds - \int_{x}^{\infty} \frac{f(s)}{(s-x)}ds \right\rbrace \, dx
\end{equation*}
\begin{equation*}
\textcolor{white}{\int_{\mathbb{R}} e^{i\beta x} \left( \frac{\partial f}{\partial |x|} \right)(x) dx} = \frac{- i \beta}{\pi}  \int_{\mathbb{R}} e^{i\beta x} \left\lbrace \int_{-\infty}^{x} \frac{f(s)}{(x-s)} ds - \int_{x}^{\infty} \frac{f(s)}{(s-x)}ds \right\rbrace \, dx
\end{equation*}
\begin{equation*}
\textcolor{white}{\int_{\mathbb{R}} e^{i\beta x} \left( \frac{\partial f}{\partial |x|} \right)(x) dx} = \frac{- i \beta}{\pi}  \int_{\mathbb{R}} f(s) \left\lbrace \int_{s}^{\infty} \frac{ e^{i\beta x}}{(x-s)} dx - \int^{s}_{-\infty} \frac{ e^{i\beta x}}{(s-x)}dx \right\rbrace
\end{equation*}
\begin{equation*}
\textcolor{white}{\int_{\mathbb{R}} e^{i\beta x} \left( \frac{\partial f}{\partial |x|} \right)(x) dx}  = \frac{-i \beta}{\pi} \left( \int_{\mathbb{R}} e^{i\beta s} f(s)ds \right) \left( \int_{0}^{\infty} \frac{ e^{i\beta y} - e^{-i\beta y}}{y} dy  \right)
\end{equation*}
\begin{equation*}
\textcolor{white}{\int_{\mathbb{R}} e^{i\beta x} \left( \frac{\partial f}{\partial |x|} \right)(x) dx}  = \frac{2 \beta}{\pi} \left( \int_{\mathbb{R}} e^{i\beta s} f(s)ds \right) \left( \int_{0}^{\infty} \frac{\sin \beta y}{y} dy  \right) .
\end{equation*}
It is well-known that
\begin{equation*}
\frac{2}{\pi}\int_{0}^{\infty} \frac{\sin \beta y}{y} dy = \frac{|\beta|}{\beta}
\end{equation*}
and thus the Fourier transform of the fractional derivative (\ref{OperatoreR}) becomes
\begin{equation*}
\int_{\mathbb{R}} e^{i\beta x} \left( \frac{\partial f}{\partial |x|} \right)(x) dx = |\beta | \int_{\mathbb{R}} e^{i\beta x} f(x) dx. 
\end{equation*}
Since 
\begin{equation}
\int_{\mathbb{R}} e^{i\beta x} \frac{t}{\pi (t^2 + x^2)} dx = e^{-t |\beta|}
\end{equation}
equation (\ref{Cdue}) immediately emerges.\\

The iterated Cauchy process defined as
\[ I_{C}(t)=C^1 (|C^2 (t)|), \quad t>0 \]
has a probability law represented by the function
\begin{equation} 
q(x,t)=\frac{2}{\pi^2} \int^{\infty}_{0} \frac{s}{s^2 + x^2} \frac{t}{t^2 + s^2} ds \quad x \in \mathbb{R}, \; t>0.
\label{DSI}
\end{equation}
We prove that
\begin{te}
The function $q=q(x,t)$ has the following explicit form
\begin{equation}
q(x,t)=\frac{2t}{\pi^2 (t^2 - x^2)} \log \frac{t}{|x|} \quad x \in \mathbb{R}, \; t>0
\label{eq:logCC}
\end{equation}
and satisfies the non-homogeneous wave equation
\begin{equation}
\frac{\partial^2 q}{\partial t^2}=\frac{\partial^2 q}{\partial x^2} - \frac{2}{\pi^2 t x^2} \quad x \in \mathbb{R}, \; t>0.
\label{eq:pdeCC}
\end{equation}
\end{te}
\paragraph{Proof:}
The integral can be easily evaluated by observing that
\[ \frac{s}{s^2+x^2} \frac{t}{t^2 + s^2} = \frac{st}{t^2 - x^2} \left\lbrace \frac{1}{x^2 + s^2} - \frac{1}{t^2+ s^2} \right\rbrace.  \]
The proof that $q$ satisfies the equation (\ref{eq:pdeCC}) can be carried out in two alternative ways.\\
By taking the second-order time derivative of \eqref{DSI} and by applying \eqref{Cuno} we have that
\begin{equation*}
\frac{\partial^2 q}{\partial t^2} = \frac{2}{\pi^2} \int_{0}^{\infty} \frac{s}{s^2 + x^2} \frac{\partial^2}{\partial t^2} \left( \frac{t}{t^2 + s^2} \right) ds
\end{equation*}
\begin{equation*}
\textcolor{white}{\frac{\partial^2 q}{\partial t^2}} = - \frac{2}{\pi^2} \int_{0}^{\infty}  \frac{s}{s^2 + x^2} \frac{\partial^2}{\partial s^2} \left( \frac{t}{t^2 + s^2} \right) ds
\end{equation*}
\begin{equation*}
\textcolor{white}{\frac{\partial^2 q}{\partial t^2}} = \frac{2}{\pi^2} \int_{0}^{\infty}  \frac{\partial}{\partial s} \left( \frac{s}{s^2 + x^2} \right) \frac{\partial}{\partial s} \left( \frac{t}{t^2 + s^2} \right) ds
\end{equation*}
\begin{equation*}
\textcolor{white}{\frac{\partial^2 q}{\partial t^2}} = \frac{2}{\pi^2} \frac{\partial}{\partial s} \left( \frac{s}{s^2 + x^2} \right) \frac{t}{t^2 + s^2} \Bigg|_{0}^{\infty} -  \frac{2}{\pi^2} \int_{0}^{\infty}  \frac{\partial^2}{\partial s^2} \left( \frac{s}{s^2 + x^2} \right) \frac{t}{t^2 + s^2} ds
\end{equation*}
\begin{equation*}
\textcolor{white}{\frac{\partial^2 q}{\partial t^2}} = -\frac{2}{\pi^2 x^2 t} + \frac{\partial^2}{\partial x^2} \frac{2}{\pi^2} \int_{0}^{\infty} \frac{s}{s^2 + x^2} \frac{t}{t^2 + s^2} ds
\end{equation*}
\begin{equation*}
\textcolor{white}{\frac{\partial^2 q}{\partial t^2}} = -\frac{2}{\pi^2 x^2 t} + \frac{\partial^2 q}{\partial x^2} .
\end{equation*}
The second proof that (\ref{eq:logCC}) satisfies (\ref{eq:pdeCC}) requires some additional attention and is based on the explicit calculations
\begin{align*}
&\frac{\partial^2 q}{\partial t^2}=\frac{2^2 t (t^2 + 3x^2)}{\pi^2 (t^2 - x^2)^3} \log \frac{t}{x} - 2 \frac{(x^2 + 3t^2)}{\pi^2t (t^2 - x^2)^2}\\
&\frac{\partial^2 q}{\partial x^2}=\frac{2^2 t (t^2 + 3x^2)}{\pi^2 (t^2 - x^2)^3} \log \frac{t}{x} + 2 \frac{t(t^2 - 5x^2)}{\pi^2 x^2 (t^2 -x^2)^2}
\end{align*}
which readily give that
\[ \frac{\partial^2 q}{\partial t^2} - \frac{\partial^2 q}{\partial x^2}=-\frac{2}{\pi^2 x^2 t}. \]
\begin{flushright}
$\blacksquare$
\end{flushright}

\begin{os}
\normalfont
The p.d.e. \eqref{eq:pdeCC} is part of a more general result of \cite[Theorem 2.1]{NA08}. 
\end{os}

The iterated Cauchy process has an equivalent representation in terms of products of independent Cauchy processes as shown in the next theorem.
\begin{te}
The following identity in distribution holds for all $t>0$
\begin{equation}
C^1(|C^2(t)|)\stackrel{i.d.}{=} \frac{1}{2} C^1(\sqrt{2t}) C^{2}(\sqrt{2t}).
\end{equation}
\end{te}
\paragraph{Proof:}
It suffices to show that for $w>0$
\begin{equation*}
Pr\left\lbrace \frac{1}{2} C^1(\sqrt{2t}) C^{2}(\sqrt{2t}) < w \right\rbrace = \frac{1}{2} + \frac{2}{\pi^2} \int_{0}^{\infty} dx \int_{0}^{\frac{2w}{x}} \frac{\sqrt{2t}}{2t + x^2} \frac{\sqrt{2t}}{2t + y^2} dy
\end{equation*}
and therefore the density $p=p(w,t)$ of the random variable $\frac{1}{2} C^1(\sqrt{2t}) C^{2}(\sqrt{2t})$ becomes
\begin{equation*}
p(w,t)=\frac{2^3 t}{2 \pi^2} \int_{0}^{\infty} \frac{x}{2t + x^2} \frac{1}{tx^2 + 2w^2} dx
\end{equation*}
\begin{equation*}
\textcolor{white}{p(w,t)} = \frac{2^2 t}{2 \pi^2 (w^2 - t^2)} \int_{0}^{\infty} \left[ \frac{x}{2t + x^2} - \frac{tx}{tx^2 + 2w^2} \right] dx
\end{equation*}
\begin{equation*}
\textcolor{white}{p(w,t)} = \frac{t}{\pi^2 (w^2 - t^2)} \log \frac{2t + x^2}{tx^2 + 2w^2} \Big|_{0}^{\infty} = \frac{2t}{\pi^2 (w^2 - t^2)} \log \frac{w}{t}.
\end{equation*}
Similar calculations must be done for $w<0$. 
\begin{flushright} 
$\blacksquare$ 
\end{flushright}

Composing more than two Cauchy processes make the research of their distribution and of the governing equation a puzzling question.\\

Two different compositions of the Brownian motion and Cauchy process are now analysed. We first obtain the p.d.e. \eqref{eq:densityCB} which is the governing equation of the process $I_{CB_H}=C(|B_H(t)|)$, $t>0$. The transition law of the process $I_{CB_H}$ is given by
\begin{equation}
p(x,t) = 2 \int_{0}^{\infty} \frac{s}{\pi (s^2 + x^2)} \frac{e^{-\frac{s^2}{2 t^{2H}}}}{\sqrt{2 \pi t^{2H}}} ds.
\label{densityProof}
\end{equation}
By performing the time derivative of \eqref{densityProof} we obtain that
\begin{equation*}
\frac{\partial p}{\partial t}= 2 \int_{0}^{\infty} \frac{s}{\pi (s^2 + x^2)} \frac{\partial}{\partial t} \left( \frac{e^{-\frac{s^2}{2 t^{2H}}}}{\sqrt{2 \pi t^{2H}}} \right) ds
\end{equation*}
\begin{equation*}
\textcolor{white}{\frac{\partial p}{\partial t}}= 2H t^{2H-1} \int_{0}^{\infty} \frac{s}{\pi (s^2 + x^2)} \frac{\partial^2}{\partial s^2} \left( \frac{e^{-\frac{s^2}{2 t^{2H}}}}{\sqrt{2 \pi t^{2H}}} \right) ds
\end{equation*}
\begin{equation*}
\textcolor{white}{\frac{\partial p}{\partial t}}= -2H t^{2H-1} \int_{0}^{\infty} \frac{\partial}{\partial s} \left( \frac{s}{\pi (s^2 + x^2)} \right) \frac{\partial}{\partial s} \left( \frac{e^{-\frac{s^2}{2 t^{2H}}}}{\sqrt{2 \pi t^{2H}}} \right) ds
\end{equation*}
\begin{equation*}
\textcolor{white}{\frac{\partial p}{\partial t}}= -2H t^{2H-1} \frac{\partial}{\partial s} \left( \frac{s}{\pi (s^2 + x^2)} \right) \frac{e^{-\frac{s^2}{2 t^{2H}}}}{\sqrt{2 \pi t^{2H}}} \Bigg|_{0}^{\infty} 
\end{equation*}
\begin{equation*}
\textcolor{white}{\frac{\partial p}{\partial t}}+ 2Ht^{2H-1} \int_{0}^{\infty} \frac{\partial^2}{\partial s^2} \left( \frac{s}{\pi (s^2 + x^2)} \right) \frac{e^{-\frac{s^2}{2 t^{2H}}}}{\sqrt{2 \pi t^{2H}}}ds.
\end{equation*}
It must be observed that the first term of the sum becomes
\begin{equation*}
-2H t^{2H-1}  \frac{x^2- s^2}{\pi (s^2 + x^2)^2} \frac{e^{-\frac{s^2}{2 t^{2H}}}}{\sqrt{2 \pi t^{2H}}} \Bigg|_{0}^{\infty} = 2Ht^{2H-1} \frac{1}{x^2}\frac{1}{\sqrt{2 \pi t^{2H}}} = \frac{2H t^{H-1}}{x^2 \pi \sqrt{2 \pi}}
\end{equation*}
and this permits us to write
\begin{equation*}
\frac{\partial p}{\partial t}= \frac{2H t^{H-1}}{x^2 \pi \sqrt{2 \pi}} - Ht^{2H-1} \frac{\partial^2}{\partial x^2}p
\end{equation*}
which concludes the proof.

\begin{os}
\normalfont
The process $I_{BC}(t)=B(|C(t)|)$, $t>0$ has probability density 
\begin{equation}
p(x,t)=2\int_{0}^{\infty} \frac{e^{-\frac{x^2}{2s}}}{\sqrt{2\pi s}} \frac{s}{\pi(s^2 + t^2)} ds, \quad x \in \mathbb{R}, \; t>0.
\end{equation}
The governing equation \eqref{eq:pdeBC} has been obtained in Nane \cite{NA08}.\\
We give here the following proof
\begin{equation*}
\frac{\partial^2 p}{\partial t^2} = 2 \int_{0}^{\infty} \frac{e^{-\frac{x^2}{2s}}}{\sqrt{2 \pi s}} \frac{\partial^2}{\partial t^2} \left( \frac{t}{\pi (s^2 + t^2)} \right) ds
\end{equation*}
\begin{equation*}
\textcolor{white}{\frac{\partial^2 p}{\partial t^2}}= -2 \int_{0}^{\infty} \frac{e^{-\frac{x^2}{2s}}}{\sqrt{2 \pi s}} \frac{\partial^2}{\partial s^2} \left( \frac{t}{\pi (s^2 + t^2)} \right) ds
\end{equation*}
\begin{equation*}
\textcolor{white}{\frac{\partial^2 p}{\partial t^2}} = 2 \int_{0}^{\infty} \frac{\partial}{\partial s} \left( \frac{e^{-\frac{x^2}{2s}}}{\sqrt{2 \pi s}} \right) \frac{\partial}{\partial s} \left( \frac{t}{\pi (s^2 + t^2)} \right) ds
\end{equation*}
\begin{equation*}
\textcolor{white}{\frac{\partial^2 p}{\partial t^2}} = 2 \frac{\partial}{\partial s} \left( \frac{e^{-\frac{x^2}{2s}}}{\sqrt{2 \pi s}} \right)  \frac{t}{\pi (s^2 + t^2)} \Bigg|_{s=0}^{s=\infty} - 2 \int_{0}^{\infty} \frac{\partial^2}{\partial s^2} \left( \frac{e^{-\frac{x^2}{2s}}}{\sqrt{2 \pi s}} \right) \frac{t}{\pi (s^2 + t^2)}ds
\end{equation*}
\begin{equation*}
\textcolor{white}{\frac{\partial^2 p}{\partial t^2}} = - \frac{1}{\pi t} \frac{\partial^2}{\partial x^2} \delta(x) - \frac{1}{2^2} \frac{\partial^4 p}{\partial x^4}.
\end{equation*}
\end{os}

\begin{os}
\normalfont
The transition law of the process $C^1(|C^2(t)|)$ can be written down in an alternative form by exploiting the subordination relationship
\begin{equation*} 
\frac{t}{\pi (x^2 + t^2)}=\int_{0}^{\infty} \frac{e^{-\frac{x^2}{2s}}}{\sqrt{2 \pi s}}\, t \,  \frac{e^{-\frac{t^2}{2s}}}{\sqrt{2 \pi s^3}}ds.
\end{equation*}
Indeed we have that
\begin{equation*}
\frac{2}{\pi^2} \int_{0}^{\infty} \frac{s}{s^2 + x^2} \frac{t}{t^2 + s^2} ds
\end{equation*}
\begin{equation*}
=2 \int_{0}^{\infty} ds \int_{0}^{\infty} dw \frac{e^{-\frac{x^2}{2w}}}{\sqrt{2 \pi w}} s \frac{e^{-\frac{s^2}{2w}}}{\sqrt{2 \pi w^3}} \int_{0}^{\infty} dz \frac{e^{-\frac{s^2}{2z}}}{\sqrt{2 \pi z}} t \frac{e^{-\frac{t^2}{2z}}}{\sqrt{2 \pi z^3}}
\end{equation*}
\begin{equation*}
= 2t \int_{0}^{\infty} dz \int_{0}^{\infty} dw \frac{e^{-\frac{x^2}{2w}}}{\sqrt{2 \pi w}} \frac{e^{-\frac{t^2}{2z}}}{\sqrt{2 \pi z^3}} \frac{1}{\sqrt{2 \pi z}} \frac{1}{\sqrt{2 \pi w^3}} \int_{0}^{\infty} ds\; s e^{-\frac{s^2}{2w}} e^{-\frac{s^2}{2z}}
\end{equation*}
\begin{equation*}
=\frac{2t}{(2\pi)^2} \int_{0}^{\infty} dz \int_{0}^{\infty} dw \frac{e^{-\frac{x^2}{2w}}}{\sqrt{w}} \frac{e^{-\frac{t^2}{2z}}}{\sqrt{z^3}} \frac{1}{\sqrt{z}} \frac{1}{\sqrt{w^3}} \frac{1}{\frac{1}{w} + \frac{1}{z}}
\end{equation*}
\begin{equation*}
=\frac{t}{2\pi^2} \int_{0}^{\infty} dz \int_{0}^{\infty} dw \frac{1}{w+z} e^{-\frac{t^2 z}{2} - \frac{x^2 w}{2}}
\end{equation*}
\begin{equation}
=\frac{1}{\pi^2} \int_{0}^{\infty} dz \int_{0}^{\infty} dw \frac{t}{zx^2 + wt^2}e^{-z-s}=\frac{1}{\pi^2} E\left\lbrace \frac{t}{x^2 Z_1 + t^2 Z_2} \right\rbrace .
\label{ExpoC}
\end{equation}
The result \eqref{ExpoC} shows that the density of the iterated Cauchy process still preserves the structure of the Cauchy distribution with exponentially distributed weights.
\end{os}

\begin{os}
\normalfont
We can give another alternative and significant representation of the distribution of the iterated Cauchy process.\\
Since
\begin{equation*}
C(\sqrt{2t}) \stackrel{i.d.}{=} \frac{1}{C(\frac{1}{\sqrt{2t}})}, \quad t>0
\end{equation*}
we can write that
\begin{equation*}
C^1(|C^2(t)|) \stackrel{i.d.}{=} \frac{1}{2\, C^1\left(\frac{1}{\sqrt{2t}}\right) C^2\left(\frac{1}{\sqrt{2t}}\right) }, \quad t>0.
\end{equation*}
Furthermore by combining the previous results we can say that
\begin{equation*}
C^1(|C^2(t)|) \stackrel{i.d.}{=} \frac{1}{C^1\left(| C^2\left(\frac{1}{t}\right)| \right)}, \quad t>0.
\end{equation*}
The last result can be also derived by observing that, for $w>0$
\begin{equation*}
Pr\left\lbrace \frac{1}{C^1\left(| C^2\left(\frac{1}{t}\right)| \right)} < w \right\rbrace = Pr\left\lbrace C^1\left( \Big| C^2\left(\frac{1}{t}\right) \Big| \right) > \frac{1}{w} \right\rbrace 
\end{equation*}
\begin{equation*}
\textcolor{white}{Pr\left\lbrace \frac{1}{C^1\left(| C^2\left(\frac{1}{t}\right)| \right)} < w \right\rbrace } = \int_{1/w}^{\infty} \int_{0}^{\infty} \frac{s}{\pi (s^2 + x^2)} \frac{\frac{1}{t}}{\pi \left( \frac{1}{t^2} + s^2 \right)} ds dx
\end{equation*}
which leads to the density of $\frac{1}{C^1\left(| C^2\left(\frac{1}{t}\right)| \right)}$ given by
\begin{equation*}
p(w,t)=\frac{2}{\pi^2} \frac{1}{w^2} \int_{0}^{\infty} \frac{s}{s^2 + \frac{1}{w^2}} \frac{\frac{1}{t}}{\frac{1}{t^2} + s^2} ds = \frac{2}{\pi^2} \int_{0}^{\infty} \frac{s}{s^2 w^2 + 1} \frac{t}{1+t^2s^2} ds
\end{equation*}
\begin{equation*}
\textcolor{white}{p(w,t)} =\frac{2}{\pi^2} \int_{0}^{\infty} \frac{\frac{z}{w}}{z^2 +1} \frac{t}{1+\frac{t^2 z^2}{w^2}} \frac{dz}{w} = \frac{2}{\pi^2} \int_{0}^{\infty} \frac{z}{z^2 + 1} \frac{t}{w^2 + t^2 z^2} dz
\end{equation*}
\begin{equation*}
\textcolor{white}{p(w,t)} = \frac{2}{\pi^2} \int_{0}^{\infty} \frac{\frac{y}{t}}{\frac{y^2}{t^2} + 1} \frac{1}{w^2 + y^2} dy = \frac{2}{\pi^2} \int_{0}^{\infty} \frac{t}{y^2 + t^2} \frac{y}{w^2 + y^2} dy.
\end{equation*}
Similar calculations hold for $w<0$.
\end{os}

\begin{table}
\centering
\begin{tabular}{c|c|c|c}
Process & Governing & Iterated & Governing \\
& equation & process & equation \\[1ex]
\hline\hline
$C(t)$ & $\frac{\partial q}{\partial t}=\frac{\partial q}{\partial |x|}$ & $C^1(|C^2(t)|)$ & $\frac{\partial^2 q}{\partial t^2}=-\frac{2}{\pi^2 t x^2} +\frac{\partial^2 q}{\partial x^2}$  \\[1ex]
 & $\frac{\partial^2 q}{\partial t^2}= - \frac{\partial^2 q}{\partial x^2}$ & & \\[1ex]
\hline
 $B_H(t)$ & $\frac{\partial p}{\partial t}=Ht^{2H-1}\frac{\partial^2 p}{\partial x^2}$ & $C(|B_H(t)|)$ & $\frac{\partial q}{\partial t}=\frac{2Ht^{2H-1}}{\pi x^2 \sqrt{2\pi}} -Ht^{2H-1}\frac{\partial^2 q}{\partial x^2}$ \\[1ex]
\hline
$B(t)$ & $\frac{\partial p}{\partial t}=\frac{1}{2}\frac{\partial^2 p}{\partial x^2}$ & $B(|C(t)|)$ & $\frac{\partial^2 q}{\partial t^2}=-\frac{1}{4}\frac{\partial^4 q}{\partial x^4} -\frac{1}{\pi t}\frac{d^2 \delta}{d x^2}$\\[1ex]
\end{tabular}
\begin{center}
\textbf{Table 2}
\end{center}
\end{table}

\paragraph{Acknowledgment}
The authors are grateful to the referee for having drawn our attention to some papers related to our work.

\end{document}